\documentclass[a4paper,12pt]{article}

\usepackage[T1]{fontenc}
\usepackage[utf8]{inputenc}
\usepackage[english]{babel}
\usepackage[inline,shortlabels]{enumitem}
\usepackage{amsmath}
\usepackage{amssymb}
\usepackage{amsthm}
\usepackage[noadjust]{cite}
\usepackage[a4paper,total={6in,9in}]{geometry}
\usepackage[colorlinks,citecolor=blue]{hyperref}
\usepackage{algpseudocode} % From algorithmicx
\usepackage[capitalize,noabbrev,nameinlink]{cleveref}
\usepackage{xcolor}
\usepackage{longtable}
\usepackage{tikz,pgfplots,pgfplotstable}
\usepackage{graphicx}
\usepackage{subcaption}
\usepackage{booktabs}

\usepgfplotslibrary{groupplots,dateplot}
\usetikzlibrary{patterns,shapes,arrows}
\pgfplotsset{compat=newest} 

\DeclareUnicodeCharacter{2212}{−}

\setlist{font=\normalfont,topsep=1ex,parsep=0ex}

\setlist[enumerate]{label=(\alph*)}

\theoremstyle{plain}
\newtheorem{proposition}{Proposition}[section]
\newtheorem{lemma}[proposition]{Lemma}
\newtheorem{theorem}[proposition]{Theorem}
\newtheorem{corollary}[proposition]{Corollary}

\theoremstyle{definition}

\newtheorem{example}[proposition]{Example}

\newtheorem{assumption}[proposition]{Assumption}
\newtheorem{algorithm}[proposition]{Algorithm}

\theoremstyle{remark}
\newtheorem{remark}[proposition]{Remark}

\newcommand{\R}{\mathbb{R}}
\newcommand{\N}{\mathbb{N}}

\newcommand{\dom}{\operatorname{dom}}
\newcommand{\weakly}{\rightharpoonup}

\newcommand{\Poi}{\operatorname{Poi}}
\def\ALM{\textbf{ALM}}
\def\ALMR{\textbf{ALMr}}
\def\SALMR{\textbf{sALMr}}

\numberwithin{equation}{section}
\numberwithin{table}{section}    
\numberwithin{figure}{section}

\crefname{figure}{Figure}{Figures}
\crefname{table}{Table}{Tables}
\Crefname{line}{Step}{Steps}

\makeatletter
\newcounter{HALG@line}
\renewcommand{\theHALG@line}{\thealgorithm.\arabic{ALG@line}}
\makeatother

\makeatletter
\newcommand{\algmargin}{\the\ALG@thistlm}
\makeatother
\algnewcommand{\parState}[1]{\State\parbox[t]{\dimexpr\linewidth-\algmargin}{\strut #1\strut}}

\let\eqref\labelcref
\usepackage{mathtools}

\newcommand\abs[1]{\left| #1\right|}
\newcommand\norm[1]{\left\Vert#1\right\Vert}
\newcommand\nnorm[1]{\Vert#1\Vert}
\newcommand\dual[2]{\left\langle#1,#2\right\rangle}
\newcommand\ddual[2]{\langle#1,#2\rangle}
\newcommand\innerp[2]{(#1,#2)}

\definecolor{todocolor}{rgb}{0.5,0.0,0.5}

\renewenvironment{proof}[1][]{\noindent\textbf{Proof{#1}: }}{\hspace*{\fill}$\Box$ \\[2mm]}

\newcommand\email[1]{\href{mailto:#1}{\texttt{#1}}}

\newcommand{\orcid}[1]{ORCID: \href{https://orcid.org/#1}{#1}}

\newcommand{\mscLink}[1]{\href{http://www.ams.org/mathscinet/msc/msc2020.html?t=#1}{#1}}

\newif\ifpreprint
\preprinttrue

\allowdisplaybreaks

\title{
		Variational Poisson Denoising via Augmented Lagrangian Methods
	}

\author{
	Christian Kanzow\thanks{%
		University of Würzburg,
		Institute of Mathematics,
		97074 Würzburg,
		Germany,
		\email{christian.kanzow@uni-wuerzburg.de},
		\orcid{0000-0003-2897-2509}
	}%
	\and
	Fabius Kr\"amer\thanks{%
		University of Regensburg,
		Faculty of Mathematics,
		93040 Regensburg,
		Germany,
		\email{Fabius.Kraemer@mathematik.uni-regensburg.de},
		\orcid{0009-0005-1422-8502}
	}%
	\and
	Patrick Mehlitz\thanks{%
		Philipps-Universität Marburg,
		Department of Mathematics and Computer Science,
		35032 Marburg,
		Germany,
		\email{mehlitz@uni-marburg.de},
		\orcid{0000-0002-9355-850X}
	} %
	\and
	Gerd Wachsmuth\thanks{%
		Brandenburgische Technische Universität Cottbus-Senftenberg,
		Institute of Mathematics,
		03046 Cottbus,
		Germany,
		\email{wachsmuth@b-tu.de},
		\orcid{0000-0002-3098-1503}
	}%
	\and
	Frank Werner\thanks{%
		University of Würzburg,
		Institute of Mathematics,
		97074 Würzburg,
		Germany,
		\email{frank.werner@uni-wuerzburg.de},
		\orcid{0000-0001-8446-3587}
	}%
}

\date{\today}

\begin{document}

\maketitle
{
	\small\textbf{\abstractname.}
		In this paper, we denoise a given noisy image by minimizing a smoothness promoting function
		over a set of local similarity measures which compare the mean of the given image and some candidate image
		on a large collection of subboxes. 
		The associated convex optimization problem possesses a huge number of constraints which are induced
		by extended real-valued functions stemming from the Kullback--Leibler divergence.
		Alternatively, these nonlinear constraints can be reformulated as affine ones,
		which makes the model seemingly more tractable.
		For the numerical treatment of both formulations of the model
		(i.e., the original one as well as the one with affine constraints),
		we propose a rather general augmented Lagrangian method
		which is capable of handling the huge amount of constraints.
		A self-contained, derivative-free, global convergence theory is provided,
		allowing an extension to other problem classes.
		For the solution of the resulting subproblems in the setting of our suggested image denoising models, 
		we make use of a suitable stochastic gradient method.
		Results of several numerical experiments are presented in order to compare
		both formulations and the associated augmented Lagrangian methods.
	\par\addvspace{\baselineskip}
}

{
	\par\addvspace{\baselineskip}
}

{
	\par\addvspace{\baselineskip}
}

\section{Introduction}\label{sec:Intro}

Denoising of images is an important task in many applications
which has received considerable attention during the last
decades, see e.g.\ \cite{FanZhangFanZhang2019,ThanhPrasathHieu2019}
for recent reviews. In this paper, we consider the problem to estimate
an image $\hat u \in L^2 \left(\Omega\right)$ on the unit square
$\Omega := (0,1)^2$ from a random set of discrete observations
$\{\omega_1,\ldots, \omega_N\} \subset \Omega$ with $N \in \N$ related to $\hat u$
as follows. We denote by
\[
Z := \sum_{i=1}^N \delta_{\omega_i},
\]
with $\delta_\omega$ being the Dirac measure centered at $\omega \in \Omega$,
the corresponding empirical process and assume that $Z$ is in fact a Poisson
point process with intensity $\hat u$, i.e., $N \in \N$ is random and
\begin{enumerate}
	\item\label{item:expected_value_of_Z} for each measurable set $A \subset \Omega$, it holds 
	$\mathbb E \left[Z(A)\right] = \int_A \hat u(\omega)\, \mathrm d \omega$, and
	\item\label{item:independence} whenever $A_1, \ldots, A_\ell \subset \Omega$ are measurable and pairwise disjoint, 
	then the random variables $Z(A_1), \ldots, Z(A_\ell)$ are stochastically independent.
\end{enumerate}
We refer to \cite{HohageWerner2016} for details on Poisson point processes and emphasize
that~\ref{item:expected_value_of_Z} and~\ref{item:independence} already imply that,
for each measurable set $A\subset\Omega$,
\begin{equation}\label{eq:poisson}
	Z(A) \sim 
	\Poi\left(  \int_A \hat u(\omega)\, \mathrm d \omega\right),
\end{equation}
i.e., the observations are essentially Poisson distributed. 
By interpreting the $L^2$-function $\hat u$ as the density of a measure $\hat U$ 
with respect to (w.r.t.) the Lebesgue measure $\mathrm d \omega$,
i.e., $\hat U (A) = \int_A \hat u(\omega) \, \mathrm d \omega$ 
for all measurable sets $A \subset \Omega$, 
the relation \eqref{eq:poisson} shows that the measure $Z$ is in fact a noisy version of $\hat U$
with $\mathbb E\left[Z(A)\right] = \hat U(A)$ for all measurable $A \subset \Omega$ 
due to~\ref{item:expected_value_of_Z}. 
In view of the Radon--Nikodým theorem, 
there is a one-to-one correspondence between $\hat u$ and 
the absolutely continuous measure $\hat U$, 
and, hence, $Z$ can be interpreted as a noisy 
(in fact Poissonian) version of $\hat u$. 
 
As the Poisson distribution is a natural model in 
applications ranging from astronomy to biophysics,
see e.g.\ \cite{VardiSheppKaufman1985,BerteroBoccacciDesderaVicidomini2009,AspelmeierEgnerMunk2015,HohageWerner2016},
the problem to estimate $\hat u$ from $Z$ has 
received considerable attention over the past decades.
Early references include \cite{SnyderWhiteHammoud1993,SnyderHelstromLantermanWhiteFaisal1995},
where explicit models for the noise occurring in 
Charge Coupled Device cameras and corresponding methods
for its removal where discussed. Since then, a major focus has been on variational
approaches, see e.g.
\cite{BonettiniRuggiero2011,BenfenatiRuggiero2013,BenvenutoCameraTheysFerrariLanteriBertero2008,
BonettiniBenfenatiRuggiero2014,FigueiredoBioucasDias2010,LeChartrandAsaki2007,
LiShenYinZhang2015,BonettiniRuggiero2012,ZanellaBoccacciZanniBertero2009}.
In all these works, a (convex) functional $J \colon L^2(\Omega) \to \overline{\R}$ of the composite form
\begin{equation}\label{eq:functional}
	\forall u\in L^2(\Omega)\colon\quad
	J(u) := G_Z(u) + f(u)
\end{equation}
is minimized, where $G_Z \colon L^2 (\Omega) \to \overline{\R}$ is a data fidelity term
measuring the discrepancy between the observations $Z$ and the
candidate image $u \in L^2 (\Omega)$, and $f \colon L^2 (\Omega) \to \overline{\R}$
is a regularization term promoting desired properties of $u$ (e.g.\ smoothness,
sparsity,$\ldots$). A natural choice for $G_Z$ is the negative log-likelihood functional
of the Poisson distribution
\begin{equation}\label{eq:nllh}
	G_Z(u) 
	:= 
	\int_\Omega u(\omega) \, \mathrm d \omega - \int_\Omega \ln(u(\omega)) \, \mathrm d Z(\omega) 
	= 
	\int_\Omega u(\omega) \, \mathrm d \omega - \sum_{i=1}^N \ln \left(u(\omega_i)\right).
\end{equation}
The previously mentioned works do all rely on this data fidelity term, and mostly differ
in the choices of $f$ and the algorithm used for minimization. 
The smoothness-promoting function $f$ 
can be chosen depending on the application. 
Famous choices include classical $L^2$-norm penalties, 
sparsity promoting penalties such as $f(u) := \sum_{i=1}^\infty |\innerp{u}{e^i}_{L^2(\Omega)}|$ 
with a complete orthonormal system or frame $\{e^i\}_{i \in \N}\subset L^2(\Omega)$, or the TV-seminorm given by 
\[
\operatorname{TV}(u) 
:= 
\sup 
\left\{ 
\int_\Omega u(\omega)\, \operatorname{div}\varphi(\omega)\, \mathrm d \omega 
\,\middle|\, 
\varphi \in C^1_c(\Omega;\R^2), \max\limits_{\omega\in\Omega}\norm{\varphi(\omega)}_2 \leq 1 
\right\}
\]
for each $u\in L^2(\Omega)$,
which equals $\operatorname{TV}(u) = \int_{\Omega} \norm{\nabla u(\omega)}_2 \mathrm d \omega$ for differentiable $u$,
see \cite{AmbrosioFuscoPallara2000} for details.
Above, $C^1_c(\Omega;\R^2)$ represents the space of all continuously differentiable functions mapping from
$\Omega$ to $\R^2$ with compact support in $\Omega$.
In the case where smooth functions are desired, also Sobolev-type penalties
\begin{equation}\label{eq:sobolev_penalty}
	f(u) := \int_{\R^2} \left(1+\|\zeta\|_2^2\right)^\theta |(\mathfrak Fu)(\zeta)|^2 \mathrm d \zeta,
\end{equation}
are suitable, where $\theta \geq 0$ and $\mathfrak F u$ denotes the Fourier transform of $u$ extended by $0$ to all of $\R^2$. 

\subsection{Problem statement}

In this paper, we will follow a different approach proposed in
\cite{FrickMarnitzMunk2012,FrickMarnitzMunk2013}. 
This is based on considering a constrained problem
\begin{equation}\label{eq:constraint}
	\min_{u\in L^2(\Omega)} \quad f(u) \quad\text{s.t.}\quad H_Z(u) \leq c
\end{equation}
for some function $H_Z\colon L^2(\Omega)\to\overline{\R}^m$ depending on $Z$ and some vector $c\in\R^m$
instead of $\min_{u\in L^2(\Omega)} J(u)$ with $J$ as in \eqref{eq:functional}. Note that, for $m = 1$ and $H_Z = \sigma G_Z$
with some $\sigma > 0$, this problem is in fact related to $\min_{u\in L^2(\Omega)} J(u)$, 
see \cite{ZanniBenfenatiBerteroRuggiero2015}. 
However, motivated by results from multiscale statistics,
we employ a different choice than (simply) the negative log-likelihood term from \eqref{eq:nllh}
for the constraint function $H_Z$ in \eqref{eq:constraint}, and aim to minimize $f$ only over
images which are everywhere locally compatible to $Z$. For the latter, we introduce a
(carefully chosen) finite system $\mathcal B \subset 2^\Omega$ of measurable regions
with positive measure in $\Omega$ (e.g., a set of square subboxes of the image), 
and consider a candidate image $u\in L^2(\Omega)$ as compatible with the data if and only if 
its mean $u_B := |B|^{-1} \int_Bu(\omega)\, \mathrm d \omega$ with the Lebesgue measure $|B|$ of $B$ 
deviates not too much from the mean $Z_B := |B|^{-1}Z(B)$ of the data $Z$ on $B$ for all $B \in \mathcal B$. 
Given the Poisson distribution of $Z_B$, deviation of $u_B$ from $Z_B$ can be made precise by means of statistical hypothesis testing, 
or as a specific instance by the local likelihood ratio test (LRT for short) statistic
\begin{equation}\label{eq:LRT_statistic}
	\forall u\in L^2(\Omega)\colon\quad
T_B \left(Z, u\right) 
:=
\sqrt{
	2 \abs{B} \left(u_B - Z_B + Z_B \ln \left(Z_B/u_B\right) \right)
}.
\end{equation}
Whenever the local LRT statistic $T_B \left(Z, u\right)$ is too large 
(which can be made precise when specifying the type 1 error of the LRT), 
the candidate image $u$ is considered incompatible with $Z$ on $B$.

This motivates the consideration of the 
\emph{variational Poisson denoising} optimization problem
\begin{equation}\label{eq:VPD}\tag{VPD}
	\min\limits_{u\in L^2(\Omega)}
	\quad
	f(u) 
	\quad
	\text{s.t.}
	\quad 
	\eta \left(Z_B, u_B\right) \leq r(\abs{B}) \quad \forall B \in \mathcal B
\end{equation}
with a function $r\colon[0,1]\to(0,\infty)$ reflecting that the right-hand side of
the constraints should -- similar to the potential number of possible regions -- 
depend on the \emph{scale} $\abs{B}$ only, and the so-called
\emph{Kullback--Leibler divergence} $\eta\colon\R^2\to\overline{\R}$ given by
\begin{equation}\label{eq:Kullback_Leibler_Divergence}
	\forall (a,b)\in\R^2\colon\quad
	\eta \left(a,b\right) 
	:=
	\begin{cases} 
		b-a + a\ln \left(a/b\right) 	& \text{if } a>0,\,b>0,\\
		b										& \text{if } a=0,\,b\geq 0,\\
		\infty									& \text{otherwise.}
	\end{cases}
\end{equation}
Note that $\eta$ is a nonnegative, convex, and lower semicontinuous function
which is continuously differentiable on $\{(a,b)\in\R^2\,|\,a,b>0\}$, 
see e.g.\ \cite{HohageWerner2016,Tsybakov2009}. 
However, $\eta$ is discontinuous precisely at the points from
$\{(a,b)\in\R^2\,|\,0\leq a\perp b\geq 0\}$ and, thus, essentially nonsmooth.
Furthermore, note the similarity to the negative log-likelihood term
$G_Z$ in \eqref{eq:nllh}, which differs from an integral over $\eta$
only by terms independent of the reconstruction candidate $u$.

The variational Poisson denoising problem is methodologically appealing from 
a statistical point of view and has a very important property, which cannot be obtained
from the variational or any other method mentioned before. 
If the function $r$ is chosen such that 
\begin{equation}\label{eq:quantile}
	\mathbb P \left[\forall B \in \mathcal B\colon\, \eta \left(Z_B, \hat u_B\right) \leq r(\abs{B})\right] \geq \alpha
\end{equation}
holds for the true image $\hat u$, i.e., if $0$ is a $ (1-\alpha)$-quantile of the random variable given by
$\sup_{B \in \mathcal B} \left[\eta \left(Z_B, \hat u_B\right) - r(\abs{B})\right]$,
see \cite{KoenigMunkWerner2020}, 
then each reconstruction $\bar u\in L^2(\Omega)$ solving \eqref{eq:VPD} satisfies automatically
\begin{equation}\label{eq:guarantee}
	\mathbb P \left[ f(\bar u) \leq f(\hat u) \right] \geq \alpha,
\end{equation}
i.e., with probability at least $\alpha$, the reconstruction $\bar u$ is at least as smooth as the true image $\hat u$. 
Further theoretical properties of a similar method in the case of Gaussian observations such as optimality for given function classes 
have been considered in \cite{GrasmairLuMunk2018,delAlamoLiMunk2021}. 
From a practical perspective, \eqref{eq:quantile} gives a clear interpretation on how to choose $r$ and makes the method
free of tuning or regularization parameters.
However, this property comes at the price that \eqref{eq:VPD} is a computationally demanding
nonsmooth convex problem, which especially suffers from a huge number of constraints. 
Exemplary, for the
setting we will consider in our numerical examples below, where $u$ is discretized by means of $266^2$
equally sized pixels, then choosing $\mathcal{B} \subset 2^{\Omega}$  as the family of all 
subsquares in the image with side length (scale) between $1$ and $64$ pixels leads already
to $3541216$ constraints, see \cref{sec:denoising_implementation}. 

In \cite{FrickMarnitzMunk2012,FrickMarnitzMunk2013}, problem \eqref{eq:VPD} has been tackled
by an ADMM approach. This requires that, in each iteration, a projection to the feasible set
\begin{align*}
	\mathcal F 
	:=&
	\left\{u \in L^2 \left(\Omega\right)\,\middle|\,\forall B\in\mathcal B\colon\,\eta \left(Z_B, u_B\right) \leq r(\abs{B})\right\}\\
	=& 
	\bigcap_{B\in \mathcal B} \left\{u \in L^2 \left(\Omega\right) \,\middle|\, \eta \left(Z_B, u_B\right) \leq r(\abs{B})\right\}
\end{align*}
has to be computed. 
In the previously mentioned works, this has been carried out by Dykstra's algorithm,
leading to a computationally poor performance. 
Let us note that if the noisy image $Z$ possesses $L^2$-regularity,
then it would be feasible to \eqref{eq:VPD}.

Further implementations of comparable problems (in the Gaussian case) have been discussed in \cite{delAlamoLiMunkWerner2020,MunkProkschLiWerner2020}, 
relying on different algorithmic frameworks such as the Chambolle--Pock algorithm, see \cite{ChambollePock2011},
or semismooth Newton methods, see \cite{HintermuellerItoKunisch2002}, 
but these approaches explicitly exploit the Gaussian and, hence, smooth structure of the corresponding optimization problem.

\subsection{Augmented Lagrangian methods}

In this paper, we will approach \eqref{eq:VPD} by means of augmented Lagrangian methods
which provide a well-established framework 
for the numerical solution of constrained optimization problems, 
see e.g.\ \cite{Bertsekas1982,BirginMartinez2014}. 
The principal idea behind those methods is to replace the minimization of a function
subject to difficult constraints by the minimization of the sum of the
associated Lagrangian function and a quadratic penalty term. Typically, the arising subproblems
are unconstrained or at least possess merely simple constraints. In contrast to penalty methods,
convergence results for augmented Lagrangian methods do not necessarily require that the
penalty parameter tends to infinity.
The augmented Lagrangian method
should be viewed as a general framework which allows an adaptation
to many different scenarios simply by taking a suitable 
and problem-dependent subproblem solver. The two standard
references mentioned above consider the situation of a 
general nonlinear program (in finite dimensions), but a 
suitable (global) convergence theory tailored for appropriate 
stationary points is also available for a couple of difficult, 
structured, and/or nonsmooth problems. This includes
situations with an abstract geometric constraint (with potentially nonconvex constraint set),
see \cite{GuoDeng2022,JiaKanzowMehlitzWachsmuth2022},
and programs with a composite objective function, see
\cite{ChenGuoLuYe2017,DeMarchiJiaKanzowMehlitz2023,DeMarchiMehlitz2024,DhingraKhongJovanovic2019,HangSarabi2021,KrugerMehlitz2022,Rockafellar2022}.
Specifically, \cite{Rockafellar2022} eliminates issues of nonsmoothness 
by exploiting smoothness properties of the Moreau
envelope in a partially convex situation. 
The fully nonsmooth setting is also discussed in
\cite{DuanXuLuZhang2019,XuYeZhang2015}, where all functions are smoothed, 
as well as in \cite{LuSunZhou2022} in the framework of so-called difference-of-convex programs.

While these references mainly deal with finite-dimensional
problems, the augmented Lagrangian approach can also be
extended to the infinite-dimensional situation. Here, we 
distinguish between the ``half'' and ``full'' infinite-dimensional
setting. Both settings share the property that the 
optimization variables belong to a Banach space,
but the former allows only finitely many inequality constraints
(possibly additional equality as well as abstract constraints),
whereas the latter allows more general functional constraints
stated in a Banach space 
(say, $ G(x) \in K $ for a mapping $ G\colon X\to Y $ between two Banach spaces $X$ and $Y$ 
as well as a closed, convex set $ K \subset Y $). 
The convergence theory for the ``half'' infinite-dimensional
setting was already considered in the seminal paper
\cite{Rockafellar1973} by Rockafellar, see also the monograph
\cite{ItoKunisch2008}. Extensions to the fully 
infinite-dimensional setting are given in 
\cite{BoergensKanzowMehlitzWachsmuth2020,BoergensKanzowSteck2019,KanzowSteckWachsmuth2018}.

We should note, however, that there exist different versions
for a realization of the augmented Lagrangian approach. 
In particular, there is the classical method with the
standard Hestenes--Powell--Rockafellar update of the 
Lagrange multipliers, and there is the safeguarded version
with a more careful updating of the multiplier estimates,
see \cite{BirginMartinez2014}. 
On the one hand, the counterexample in 
\cite{KanzowSteck2017} shows that there cannot exist a 
satisfactory global convergence theory for the classical 
method, at least not in the nonconvex setting, while the
existing convergence theory for the safeguarded version is
rather complete in the sense that it has all desirable 
(and realistic) properties. On the other hand, for
convex problems, there exists a convergence theory for 
the classical approach even with a constant penalty parameter.
This result was established by Rockafellar in \cite{Rockafellar1973},
even for the ``half'' infinite-dimensional
setting, and is based on the duality of the augmented Lagrangian
and the proximal point method, see \cite{Rockafellar2022} as well. 
In particular, this convergence
theory is based on the existence of Lagrange multipliers.

\subsection{Our contributions}

As mentioned earlier, the purpose of this article is to study the numerical solution
of \eqref{eq:VPD} with the aid of augmented Lagrangian methods.
Let us note that \eqref{eq:VPD} is covered by the ``half'' infinite-dimensional setting, and
due to the huge number of constraints in \eqref{eq:VPD}, the augmentation approach seems
to be perfectly suitable to tackle the problem computationally. 
In the course of the paper, we also suggest a reformulation of the constraints in \eqref{eq:VPD}
as purely affine inequalities, and the latter is covered by the ``half'' infinite-dimensional
setting as well. 
Though both models are convex, we provide a purely primal convergence
theory for a safeguarded augmented Lagrangian method 
in a more general nonconvex setting. We assume, however, that we
are able to find an approximate global minimum of the resulting subproblems. 
This is a realistic scenario for the convex optimization problem
\eqref{eq:VPD} and its aforementioned reformulation 
but might also be applicable in some other situations (e.g., think
of disjunctive constraint systems composed of finitely many convex branches). 
In contrast to the existing literature, our analysis is derivative-free and
works under minimal continuity assumptions on the data.
This is rather important as the constraints in \eqref{eq:VPD} are modeled
with the aid of the discontinuous function $\eta$ from \eqref{eq:Kullback_Leibler_Divergence}.
Moreover, we stress that our convergence theory is independent 
of any assumption regarding the existence (or uniqueness or boundedness) 
of Lagrange multipliers.
Finally, let us mention that the method has favorable convergence properties 
even in the case where the constraints are inconsistent 
since it still provides limits which minimize a suitable measure 
for the constraint violation. 

The suggested safeguarded augmented Lagrangian method is then applied to the numerical
solution of \eqref{eq:VPD} as well as its reformulation with affine constraints.
As the latter is guaranteed to possess Lagrange
multipliers at its minimizers, 
we also apply to it the classical augmented Lagrangian method
(i.e., we abstain from safeguarding) with a constant penalty parameter which is reasonable
due to the analysis in \cite{Rockafellar1973}. We compare these three numerical approaches
by means of several different test instances and document the results.

The paper is organized in the following way.
In the remainder of this introductory section,
we comment on the notation which will be used throughout the paper.
The safeguarded augmented Lagrangian
method of our interest is stated and analyzed in \cref{sec:ALM}, where we
consider nonsmooth problems with finitely many
inequality constraints, a general operator equation
(representing, e.g., a partial differential equation),
as well as an abstract constraint set such that the
associated augmented Lagrangian subproblems can be solved
up to approximate global optimality.
\cref{sec:experiments} is dedicated to the numerical solution of the Poisson denoising model
with the aid of augmented Lagrangian methods.
In \cref{sec:reformulation}, we derive a reformulation of \eqref{eq:VPD}
which merely possesses (finitely many) affine inequality constraints.
The implementation of our experiments, the exploited subproblem solver, and the way
we are documenting our results are 
discussed in \cref{sec:denoising_implementation,sec:denoising_SGD,sec:documentation}, respectively,
before the numerical performance of augmented Lagrangian methods when applied to \eqref{eq:VPD}
and its reformulation is illustrated and compared based on several different test instances
in \cref{sec:denoising_numerical_results}.
We conclude with some final remarks in \cref{sec:conclusions}.

\subsection{Notation}

Let $\R$ and $\R_+$ denote the sets of all real numbers and nonnegative real numbers, respectively.
We make use of $\overline{\R}:=\R\cup\{\infty\}$.
Throughout the paper, for a given finite set $D$, $\#D$ is used to denote the cardinality of $D$.
Let $n\in\N$ be a positive integer.
For vectors $x,y\in\R^n$, $\max(x,y)\in\R^n$ denotes the componentwise maximum of $x$ and $y$.
For any $p\in[1,\infty]$, the $\ell_p$-norm of $x\in\R^n$ will be denoted by $\norm{x}_p$.

Whenever $X$ is a Banach space, its norm will be denoted by $\norm{\cdot}_X\colon X\to [0,\infty)$
if not stated otherwise. Strong and weak convergence of a sequence $\{x^k\}\subset X$ to $x\in X$
are represented by $x^k\to x$ and $x^k\weakly x$, respectively.
If $K\subset\N$ is a set of infinite
cardinality, we make use of $x^k\to_K x$ ($x^k\weakly_K x$) in order to express that the
subsequence $\{x^k\}_{k\in K}$ converges (converges weakly) to $x$ as $k$ tends to $\infty$
in $K$ (which we denote by $k\to_K\infty$ for brevity).
The (topological) dual space of $X$ will be represented by $X^*$, and the associated dual
pairing is then denoted by $\dual{\cdot}{\cdot}_X\colon X^*\times X\to\R$.
Let $Y$ be another Banach space. 
If $h\colon X\to Y$ is Fr\'{e}chet differentiable at $x\in X$, $h'(x)\colon X\to Y$ denotes the derivative
of $h$ at $x$. Similarly, if $X_1$ and $X_2$ are Banach spaces such that $X=X_1\times X_2$,
and if $h$ is Fr\'{e}chet differentiable at $x:=(x_1,x_2)\in X$, $h'_{x_1}(x)\colon X_1\to Y$ denotes the
partial derivative w.r.t.\ $x_1$ of $h$ at $x$.
The inner product in a Hilbert space $H$ will be represented by 
$\innerp{\cdot}{\cdot}_H\colon H\times H\to\R$.

For an arbitrary function $\varphi\colon X\to\overline\R$ defined on a Banach space $X$, 
the set $\dom\varphi:=\{x\in X\,|\,\varphi(x)<\infty\}$ is referred to as the domain of $\varphi$. 
Whenever $\varphi$ is convex and $\bar x\in\dom\varphi$ is chosen arbitrarily, the set
\[
	\partial\varphi(\bar x)
	:=
	\{
		\xi\in X^*
		\,|\,
		\forall x\in\dom\varphi\colon\,
		\varphi(x)\geq \varphi(\bar x)+\dual{\xi}{x-\bar x}_X
	\}
\]
is called the subdifferential (in the sense of convex analysis) of $\varphi$ at $\bar x$.

For an integer $d\in\N$, a bounded open set $\Omega\subset\R^d$, and $p\in[1,\infty)$, 
$L^p(\Omega)$ denotes the Lebesgue space of (equivalence classes of)
measurable functions $u\colon\Omega\to\R$ such that $\Omega\ni\omega\mapsto\abs{u(\omega)}^p\in\R$
is integrable, and is equipped with the standard norm which we denote by 
$\norm{\cdot}_p\colon L^p(\Omega)\to [0,\infty)$.
Note that it will be clear from the context 
whether $\norm{\cdot}_p$ is taken in $\R^n$ or $L^p(\Omega)$.

\section{An augmented Lagrangian method for nonsmooth optimization problems}\label{sec:ALM}

In this section, we address the algorithmic treatment of the rather abstract optimization problem
\begin{equation}\label{eq:nonsmooth_problem}\tag{P}
	\min\limits_{x\in X} \ f(x) \quad \text{s.t.} \quad g(x) \leq 0, \ h(x) = 0, \ x\in C,
\end{equation}
where $f\colon X\to\overline{\R}$, $g\colon X\to\overline{\R}^m$, and
$h\colon X\to Y$ are given functions and $C\subset X$ is a weakly sequentially closed set.
Moreover, $X$ is a reflexive Banach space and $Y$ is a Hilbert space,
which we identify with its dual, i.e., $Y\cong Y^*$.
Throughout this section, 
let $\mathcal F\subset X$ denote the feasible set of \eqref{eq:nonsmooth_problem}.
For later use, we introduce $\dom g:=\bigcap_{i=1}^m\dom g_i$.
Here, $g_1,\ldots,g_m$ are the component functions of $g$.
As a minimal requirement, we need 
\begin{equation}\label{eq:not_fully_trivial}
	\dom f\cap\dom g\cap C\neq\emptyset,
\end{equation}
which will be assumed to be a standing assumption in this section.
Observe that \eqref{eq:not_fully_trivial} does not necessarily rule out 
$\dom f\cap\mathcal F=\emptyset$.
However, $\dom f\cap\mathcal F\neq\emptyset$ 
clearly implies validity of \eqref{eq:not_fully_trivial}.

In contrast to the standard setting of nonlinear programming, we abstain from
demanding any differentiability properties of the data functions.
However, we assume that the functions
$f,g_1,\ldots,g_m\colon X\to\overline{\R}$ are weakly sequentially lower semicontinuous, while the function
$h$ is weakly-strongly sequentially continuous in the sense that
\[
	\forall \{x^k\}\subset X,\,\forall x\in X\colon\quad
	x^k\weakly  x \quad \text{in} \; X
	\quad\Longrightarrow\quad
	h(x^k)\to h(x)	\quad \text{in} \; Y.
\]
Note that at least continuity of the function $h$ is indispensable in order to
guarantee that $\mathcal F$ is closed. 
The assumptions from above already guarantee that $\mathcal F$ is weakly sequentially closed.
Together with $\dom f\cap\mathcal F\neq\emptyset$ 
and the weak sequential lower semicontinuity of the objective functional $f$, 
this can be interpreted as a minimal requirement in constrained optimization in order to ensure 
that the underlying optimization problem \eqref{eq:nonsmooth_problem} possesses a solution.
This would be inherent whenever $\mathcal F$ is, additionally, bounded or $f$ is, 
additionally, coercive as standard arguments show.

\subsection{A chain rule for lower semicontinuity}

Before we can start with the presentation of the augmented Lagrangian method and its convergence analysis,
we need to prepare conditions which guarantee that the composition of a
(weakly sequentially) lower semicontinuous function and a continuous function is 
(weakly sequentially) lower semicontinuous again.
Such a criterion is presented in the upcoming lemma.

\begin{lemma}\label{lem:chain_rule_lower_semicontinuity}
	For some Banach space $X$, let $\varphi\colon X\to\overline{\R}$ be weakly sequentially lower semicontinuous
	and let $\psi\colon\R\to\R$ be a continuous and monotonically increasing
	function.
	Then $\psi\circ\varphi\colon X\to\overline{\R}$ defined via
	\[
		\forall x\in X\colon\quad
		(\psi\circ\varphi)(x)
		:=
		\begin{cases}
			\psi(\varphi(x))				&\text{if } \varphi(x)<\infty,\\
			\lim_{t\to\infty}\psi(t)		&\text{if } \varphi(x)=\infty
		\end{cases}
	\]
	is weakly sequentially lower semicontinuous.
\end{lemma}
\begin{proof}
	Choose $\{x^k\}\subset X$ and $\bar x\in X$ with $x^k\weakly\bar x$ arbitrarily and
	pick an infinite set $K\subset\N$
	such that 
	\[
		\alpha:=\liminf_{k\to\infty}(\psi\circ\varphi)(x^k)=\lim_{k\to_K\infty}(\psi\circ\varphi)(x^k).
	\]
	First, we argue that $\alpha$ cannot attain the value $-\infty$.
	Indeed,
	weak sequential lower semicontinuity of $\varphi$ yields the estimate
	$-\infty < \varphi(\bar x) \le \liminf_{k\to\infty} \varphi(x^k)$,
	so we infer that $\{\varphi(x^k)\}$ is bounded from below.
	Consequently,
	$\{(\psi\circ\varphi)(x^k)\}$
	is also bounded from below and, thus, $\alpha > -\infty$.
	In the case $\alpha=\infty$, we automatically have $\alpha\geq(\psi\circ\varphi)(\bar x)$ and, thus,
	there is nothing to show.
	Thus, we assume $\alpha\in\R$. By weak sequential lower semicontinuity of $\varphi$, we have
	$\beta:=\liminf_{k\to_K\infty}\varphi(x^k)\geq\varphi(\bar x)$. Pick an infinite set
	$K'\subset K$ such that $\lim_{k\to_{K'}\infty}\varphi(x^k)=\beta$. In the case where
	$\beta\in\R$ holds, $\varphi(\bar x)$ and the tail of the sequence $\{\varphi(x^k)\}_{k\in K'}$ are finite, so we find 
	\begin{align*}
		\alpha
		=
		\lim\limits_{k\to_{K'}\infty}(\psi\circ\varphi)(x^k)
		=
		\lim\limits_{k\to_{K'}\infty}\psi(\varphi(x^k))
		=
		\psi(\beta)
		\geq
		\psi(\varphi(\bar x))
		=
		(\psi\circ\varphi)(\bar x)
	\end{align*}
	by continuity and monotonicity of $\psi$.	
	Next, suppose that $\beta=\infty$ holds. 
	In the case where $\{\varphi(x^k)\}_{k\in K'}$ equals $\infty$ along the tail of the
	sequence, we find
	\[
		\alpha
		=
		\lim\limits_{k\to_{K'}\infty}(\psi\circ\varphi)(x^k)
		=
		\lim\limits_{t\to\infty}\psi(t)
		\geq
		(\psi\circ\varphi)(\bar x)
	\]
	by monotonicity of $\psi$.
	Otherwise, there is an infinite set $K''\subset K'$ such that we have $\{\varphi(x^k)\}_{k\in K''}\subset\R$.
	However, $\beta=\infty$ yields $\lim_{k\to_{K''}\infty}\varphi(x^k)=\infty$.
	Hence, by definition of the composition, we find
	\[
		\alpha
		=
		\lim\limits_{k\to_{K''}\infty}(\psi\circ\varphi)(x^k)
		=
		\lim\limits_{k\to_{K''}\infty}\psi(\varphi(x^k))
		=
		\lim\limits_{t\to\infty}\psi(t)
		\geq
		(\psi\circ\varphi)(\bar x)
	\]
	by continuity and monotonicity of $\psi$.
	This completes the proof.
\end{proof}

We would like to note that, in general, for a 
(weakly sequentially) lower semicontinuous function
$\varphi\colon X\to\overline{\R}$, the mappings $x\mapsto|\varphi(x)|$ 
and $x\mapsto \varphi^2(x)$ are not (weakly sequentially)
lower semicontinuous (exemplary, choose $X:=\R$ and set $\varphi(x):=-1$ for
all $x\leq 0$ and $\varphi(x):=0$ for all $x>0$). 
Observe that the absolute value function and the square are not monotonically increasing, i.e.,
the assumptions of \cref{lem:chain_rule_lower_semicontinuity} are not satisfied
in this situation.

We comment on a typical setting where \cref{lem:chain_rule_lower_semicontinuity}
applies.
\begin{example}\label{ex:lower_semicontinuity_of_ALM}
	For each $\alpha>0$ and $\beta\in\R$, the function $\psi\colon\R\to\R$
	given by $\psi(t):={\max}^2(0,\alpha t+\beta)$ for each $t\in\R$ is
	continuous, monotonically increasing, and satisfies $\lim_{t\to\infty}\psi(t)=\infty$.
	Thus, for each Banach space $X$ 
	and each weakly sequentially lower semicontinuous function $\varphi\colon X\to\overline\R$, the
	composition $\psi\circ\varphi\colon X\to\overline{\R}$ given by
	\[
		\forall x\in X\colon\quad
		(\psi\circ\varphi)(x)
		:=
		\begin{cases}
			\psi(\varphi(x))	&	\text{if } \varphi(x)<\infty,\\
			\infty				&	\text{if } \varphi(x)=\infty
		\end{cases}
	\]
	is weakly sequentially lower semicontinuous as well by \cref{lem:chain_rule_lower_semicontinuity}.
	
	We also note that this particular function $\psi$ is convex.
	Thus, keeping the monotonicity of $\psi$ in mind,
	whenever $\varphi$ is convex, 
	then the composition $\psi\circ\varphi$ is convex as well.
\end{example}

\subsection{Statement of the algorithm}\label{sec:ALM_algorithm}

Let $L\colon X\times\R^m_+\times Y\to \overline{\R}$ denote the 
Lagrangian function associated with \eqref{eq:nonsmooth_problem} which is given by
\[
	L(x,\lambda,\mu):=f(x)+\lambda^\top g(x)+\innerp{\mu}{h(x)}_Y
\]
for $x\in X$, $\lambda\in\R^m_+$, and $\mu\in Y$.
For the construction of our solution method, we make use of the corresponding augmented
Lagrangian function $L_\rho\colon X\times\R^m_+\times Y\to\overline{\R}$ 
associated with \eqref{eq:nonsmooth_problem} which is given by
\begin{equation}\label{eq:augmented_lagrangian}
	\begin{aligned}
   		L_{\rho} (x, \lambda, \mu) 
  		 := 
  		 f(x) 
  		 &+ 
 		  \frac{1}{2 \rho} \sum_{i=1}^m 
 		  	\left( {\max}^2 \left(0, \lambda_i + \rho g_i(x) \right) - \lambda_i^2 \right)
 		  \\
 		  &+
 		  \innerp{\mu}{h(x)}_Y + \frac{\rho}{2}\norm{ h(x) }^2_Y
 	\end{aligned}
\end{equation}
for all
$x\in X$, $\lambda\in\R^m_+$, and $\mu\in Y$,
where $\rho>0$ is a given penalty parameter.
Roughly speaking, the latter results from $L$ by adding a standard quadratic penalty
for the constraints $g(x)\leq 0$ and $h(x)=0$ appearing in \eqref{eq:nonsmooth_problem}.
Within our algorithmic framework, the function $L_\rho$ has to be minimized
w.r.t.\ $x$, which means that the term $-\tfrac{1}{2\rho}\norm{\lambda}^2_2$ could be removed
from the definition of $L_\rho$. However, for some of the proofs we are going to provide, it will
be beneficial to keep this shift. 
We would like to point the reader's attention to the fact that the function $L_\rho(\cdot,\lambda,\mu)$ is
weakly sequentially lower semicontinuous for each $\lambda\in\R^m_+$ and $\mu\in Y$ due to 
\cref{lem:chain_rule_lower_semicontinuity}, \cref{ex:lower_semicontinuity_of_ALM}, and the
fact that the function $h$ is weakly-strongly sequentially continuous.
Furthermore, we would like to point out that the abstract constraint set $C$ is not 
incorporated into the definitions of $L$ and $L_\rho$ on purpose.

\begin{remark}\label{rem:convex_setting}
Whenever \eqref{eq:nonsmooth_problem} is a convex optimization
problem, i.e., whenever the functions $f,g_1,\ldots,g_m$ are convex while $h$ is affine, then,
for each $\lambda\in\R^m_+$ and $\mu\in Y$, $L_\rho(\cdot,\lambda,\mu)$ is a convex function as well 
by monotonicity and convexity of $t\mapsto{\max}^2(0,\alpha t+\beta)$ for each $\alpha>0$ and $\beta\in\R$.
\end{remark}

For some penalty parameter $\rho>0$, we introduce a function $V_\rho\colon X\times\R^m_+\to\overline\R$
by means of
\[
	V_\rho(x,\lambda)
	:=
	\begin{cases}
		\max\bigl(\norm{\max(g(x),-\lambda/\rho)}_\infty,\norm{h(x)}_Y\bigr)
			&	\text{if } x \in \dom g,\\
		\infty	
			&	\text{if } x\notin\dom g
	\end{cases}
\]
for all $x\in X$ and $\lambda\in\R^m_+$. 
Right from the definition of $V_\rho$, we obtain
\[
	V_\rho(x,\lambda)=0
	\quad \Longleftrightarrow \quad
	g(x)\leq 0,\,\lambda\geq 0,\,\lambda^\top g(x)=0,\,h(x)=0 ,
\]
i.e., $V_\rho$ can be used to measure feasibility of $x$ for \eqref{eq:nonsmooth_problem}
w.r.t.\ the constraints induced by $g$ and $h$ as well as validity of the complementarity-slackness 
condition w.r.t.\ the inequality constraints.

In \cref{alg:ALM}, we state a pseudo-code which describes our method.

\begin{algorithm}[Safeguarded Augmented Lagrangian Method for \eqref{eq:nonsmooth_problem}]\leavevmode
	\label{alg:ALM}
	\begin{algorithmic}[1]
		\Require bounded sets $B_m\subset\R^m_+$ and $B_Y\subset Y$,
 			starting point $ (x^0, \lambda^0, \mu^0) \in C\times \R^m_+ \times Y $, 
      		initial penalty parameter $ \rho_0 > 0 $, 
     		 parameters $ \tau \in (0, 1)$,  $\gamma > 1 $
		\parState{Set $k := 0$.}
		\While{a suitable termination criterion is violated at iteration $ k $}\label{item:termination_ALM} 
		\parState{Choose $v^k\in B_m$ and $w^k\in B_Y$.
			}\label{item:choice_of_approximate_multipliers}
		\parState{Compute $ x^{k + 1} \in C$ 
 			as an approximate solution of the optimization problem
      		\begin{equation}\label{eq:ALM_subproblem}
         		\min_{x\in X} \quad L_{\rho_k} (x, v^k, w^k) \quad \text{s.t.} \quad  x \in C.
      		\end{equation}}\label{item:solve:subproblem_ALM}
		\parState{Set
      		\begin{equation}\label{eq:update_multiplier}
         			\lambda^{k + 1} 
         			:=  \max \bigl( 0, v^k + \rho_k g(x^{k + 1}) \bigr),
         			\qquad
         			\mu^{k + 1}
         			:=  w^k + \rho_k h (x^{k + 1}). 
      		\end{equation}}\label{item:adjust_multiplier_ALM}
		\parState{If $ k = 0 $ or the condition
      		\begin{equation}\label{eq:update_penalty_ALM}
      			V_{\rho_k}(x^{k+1},v^k)
      			\leq
      			\tau\,V_{\rho_{k-1}}(x^k,v^{k-1})
     		\end{equation}
     		 holds, set $ \rho_{k + 1} := \rho_k $, otherwise set
      		$ \rho_{k + 1} := \gamma \rho_k $.
      		}\label{item:update_penalty_ALM} 
      	\parState{Set $ k \leftarrow k + 1 $.}
		\EndWhile
		\State\Return $x^k$
	\end{algorithmic}
\end{algorithm}

In \cref{alg:ALM}, the quantities $v^k$ and $w^k$ play the role of Lagrange
multiplier estimates. By construction, the sequences $\{v^k\}$ and $\{w^k\}$
remain bounded throughout a run of the algorithm while this does not necessarily hold 
true for $\{\lambda^k\}$ and $\{\mu^k\}$. Note that the classical augmented
Lagrangian method could be recovered from \cref{alg:ALM} by replacing
$v^k$ and $w^k$ by $\lambda^k$ and $\mu^k$ everywhere, respectively. 
However, the so-called safeguarded variant from \cref{alg:ALM} has 
been shown to possess better global convergence properties than the classical method,
see, e.g., \cite{KanzowSteck2017} for details.
Typically, $\{v^k\}$ and $\{w^k\}$ are iteratively constructed during the run of
\cref{alg:ALM}. Exemplary, one can choose $B_m$ as the (very large) box $[0,\mathbf{v}]$
for some $\mathbf v\in\R^m$ satisfying $\mathbf v>0$ and define $v^k$ as the projection of $\lambda^k$ 
onto this box in \cref{item:choice_of_approximate_multipliers}. 
A similar procedure is possible for the choice of $w^k$.
This way, \cref{alg:ALM} is likely
to parallel the classical augmented Lagrangian method if the sequences $\{\lambda^k\}$ and
$\{\mu^k\}$ remain bounded.

Assuming for a moment that all involved data functions are smooth, the derivative w.r.t.\
$x$ of $L_\rho$ from \eqref{eq:augmented_lagrangian} is given by
\begin{align*}
	(L_\rho)'_x(x,\lambda,\mu)
	=
	f'(x)
	+
	\sum\limits_{i=1}^m\max(0,\lambda_i+\rho\,g_i(x))\, g'_i(x)
	+
	h'(x)^*[\mu+\rho h(x)].
\end{align*}
Thus, the updating rule for the multipliers in \eqref{eq:update_multiplier} yields
\begin{equation}\label{eq:motivation_update_multipliers}
	(L_{\rho_k})'_x(x^{k+1},v^k,w^k)=L'_x(x^{k+1},\lambda^{k+1},\mu^{k+1}),
\end{equation}
which is the basic idea behind \cref{item:adjust_multiplier_ALM}. 
Note that a similar formula as \eqref{eq:motivation_update_multipliers} can be obtained
in terms of several well-known concepts of subdifferentiation whenever a suitable chain
rule applies.

Finally, let us mention that in \cref{item:update_penalty_ALM}, the penalty parameter is
increased whenever the new iterate $(x^{k+1},v^k,w^k)$ 
is not (sufficiently) better from the viewpoint of feasibility (and complementarity) than the 
old iterate $(x^k,v^{k-1},w^{k-1})$. Note that our choice for the infinity norm
in the definition of $V_\rho$ is a matter of taste since all norms are equivalent
in finite-dimensional spaces. 
However, this particular error measure $V_\rho$ keeps track of the largest violation
of the feasibility and complementarity condition w.r.t.\ \emph{all} inequality constraints,
which is why we favor it here.

For further information about (safeguarded) augmented Lagrangian methods in nonlinear
programming, we refer the interested reader to the monograph \cite{BirginMartinez2014}.

\subsection{Convergence to global minimizers}\label{sec:ALM_convergence}

In this subsection, we provide a convergence analysis for \cref{alg:ALM} where we assume that 
in \cref{item:solve:subproblem_ALM}, the subproblem \eqref{eq:ALM_subproblem} is solved up to (approximate) global optimality.
Exemplary, this is possible whenever \eqref{eq:nonsmooth_problem} is a convex program, see \cref{rem:convex_setting},
but also in more general situations where \eqref{eq:nonsmooth_problem} is of special structure, e.g.\ if the feasible set
can be decomposed into a moderate number of convex branches while the objective function is convex.
Within the assumption below, which will be standing throughout
this section, we quantify the requirements regarding the subproblem solver.

\begin{assumption}\label{ass:approximate_solution}
	In each iteration $k\in\N$ of \cref{alg:ALM}, the approximate solution $x^{k+1}\in C$ 
	of \eqref{eq:ALM_subproblem} satisfies
	\begin{equation}\label{eq:eps_optimality}
		\forall x\in C\colon\quad
		L_{\rho_k}(x^{k+1},v^k,w^k)-\varepsilon_k
		\leq
		L_{\rho_k}(x,v^k,w^k)
	\end{equation}
	where $\varepsilon_k\geq 0$ is some given constant.
\end{assumption}

Typically, the inexactness parameter $\varepsilon_k$ in \cref{ass:approximate_solution} is
chosen to be positive. While $\varepsilon_k:=0$ corresponds to the situation where the
subproblems \eqref{eq:ALM_subproblem} are solved exactly, 
we will see that the augmented Lagrangian technique generally
works fine if only approximate solutions of the subproblems are computed. This also has the
advantage that whenever $\inf_x\{L_{\rho_k}(x,v^k,w^k)\,|\,x\in C\}$ is finite, 
then one can always find
points $x^{k+1}$ satisfying \eqref{eq:eps_optimality} for arbitrarily small $\varepsilon_k>0$
while an exact global minimizer may not exist.
Furthermore, we note that, due to \eqref{eq:not_fully_trivial},
$L_{\rho_k}(x^{k+1},v^k,w^k)<\infty$ holds for each $k\in\N$, i.e.,
$x^{k+1}\in \dom f\cap\dom g\cap C$ is valid for each computed iterate.
Finally, it is worth mentioning that validity of \eqref{eq:eps_optimality} 
guarantees that $L_{\rho_k}(\cdot,v^k,w^k)$ is bounded from below on $C$.

Throughout the section, we make use of the following lemma.

\begin{lemma}\label{lem:upper_bound_AL_values}
	Let $v\in\R^m$, $w\in Y$, and $\rho>0$ as well as a feasible point $x\in \mathcal F$ of
	\eqref{eq:nonsmooth_problem} be arbitrary. 
	Then $L_\rho(x,v,w)\leq f(x)$ is valid.
\end{lemma}
\begin{proof}
	Due to $h(x)=0$ and by definition of the augmented Lagrangian function $L_\rho$ from
	\eqref{eq:augmented_lagrangian}, we find
	\[
		L_\rho(x,v,w)
		=
		f(x)+\frac1{2\rho}\sum\limits_{i=1}^m\left({\max}^2(0,v_i+\rho\,g_i(x))-v_i^2\right),
	\]
	i.e., in order to show the claim, it is sufficient to verify
	${\max}^2(0,v_i+\rho\,g_i(x))\leq v_i^2$ for all $i\in\{1,\ldots,m\}$.
	Thus, fix $i\in\{1,\ldots,m\}$ arbitrarily. 
	In the case $v_i+\rho\,g_i(x)\leq 0$, we find ${\max}^2(0,v_i+\rho\,g_i(x))=0\leq v_i^2$.
	Conversely, $v_i+\rho\,g_i(x)>0$ yields
	$0\leq v_i+\rho\,g_i(x)\leq v_i$ since $g_i(x)\leq 0$ is valid by feasibility of $x$ for
	\eqref{eq:nonsmooth_problem}, so by monotonicity of the square on the nonnegative real line,
	${\max}^2(0,v_i+\rho\,g_i(x))\leq v_i^2$ follows.
\end{proof}

Let us now start with the convergence analysis associated with \cref{alg:ALM}.
Therefore, we first study issues related to the feasibility of accumulation points.

\begin{proposition}\label{prop:feasibility_of_accumulation_points}
	Let \eqref{eq:not_fully_trivial} hold. 
	Assume that \cref{alg:ALM} produces a sequence $\{x^k\}$ such that \cref{ass:approximate_solution}
	holds for some bounded sequence $\{\varepsilon_k\}$, 
	and let $\{\rho_k\}$ and $\{v^k\}$ be the
	associated sequences of penalty parameters and Lagrange multiplier estimates associated with
	the inequality constraints in \eqref{eq:nonsmooth_problem}, respectively.
	Let the subsequence $\{x^{k+1}\}_{k\in K}$ and $\bar x\in X$ be chosen such that $x^{k+1}\weakly_K\bar x$.
	Then $\bar x$ is a global minimizer of the optimization problem
	\begin{equation}\label{eq:feasibility_problem}
		\min\limits_{x\in X} \quad \frac12\norm{\max(g(x),0)}_2^2+\frac12\norm{h(x)}_Y^2 
		\quad \textup{s.t.} \quad 
		x\in \dom f\cap C,
	\end{equation}
	where, for each $x\notin\dom g$, $\frac12\norm{\max(g(x),0)}_2^2:=\infty$.
	Furthermore, whenever $\mathcal F\neq\emptyset$,
	$\bar x\in\dom f\cap\mathcal F$ holds,
	and we have $V_{\rho_k}(x^{k+1},v^k)\to_K 0$.
\end{proposition}

\begin{proof}
	Let us start with the observation that, due to $x^{k+1}\in \dom f\cap C$ for each $k\in K$,
	weak sequential lower semicontinuity of $f$ and weak sequential closedness of $C$
	guarantee $\bar x\in\dom f\cap C$, 
	i.e., $\bar x$ is feasible to \eqref{eq:feasibility_problem}.
	
	To proceed, we distinguish two cases.
	
	\emph{Case 1}: Suppose that $\{\rho_k\}$ remains bounded. Then \cref{item:update_penalty_ALM}
	yields that $\rho_k$ remains constant on the tail of the sequence, i.e., there is some
	$k_0\in\N$ such that $\rho_k=\rho_{k_0}$ is valid for all $k\in\N$ satisfying $k\geq k_0$.
	Particularly, condition \eqref{eq:update_penalty_ALM} is satisfied for all $k\geq k_0$,
	which immediately yields $V_{\rho_k}(x^{k+1},v^k)\to 0$ due to $\{x^{k+1}\}\subset\dom g$.
	On the one hand, we infer $h(x^{k+1})\to 0$ and, by weak-strong sequential continuity of $h$,
	$h(x^{k+1})\to_K h(\bar x)$ on the other hand. 
	By uniqueness of the limit, $h(\bar x)=0$ follows.
	Due to boundedness of $\{v^k\}$, we may also assume w.l.o.g.\ 
	that $v^k\to_K\bar v$ is valid for some $\bar v\in\R^m$.
	The componentwise weak sequential lower semicontinuity of $g$ yields 
	$\max(g(\bar x),-\bar v/\rho_{k_0})\leq 0$ in the light of 
	\eqref{eq:update_penalty_ALM}, i.e., $g(\bar x)\leq 0$ follows. 
	Hence, $\bar x\in\mathcal F$ has been shown, 
	i.e., $\bar x$ is feasible to \eqref{eq:nonsmooth_problem}. 
	The objective value of \eqref{eq:feasibility_problem}
	associated with $\bar x$ is $0$, 
	so that $\bar x$ is a global minimizer of this problem as well.
	
	\emph{Case 2}: Now, assume that $\{\rho_k\}$ is not bounded.
	Then, by construction, we have $\rho_k\to\infty$.
	
	Fix an arbitrary point $x\in \dom f\cap\dom g\cap C$
	(such a point exists due to our standing assumption \eqref{eq:not_fully_trivial}).
	Observe that \cref{ass:approximate_solution} 
	and the definition of the augmented Lagrangian function give
	\begin{equation}\label{eq:some_interesting_estimate}
		\begin{aligned}
		&f(x^{k+1}) 
   		+ 
   		\frac{1}{2 \rho_k} \sum_{i=1}^m 
   			\left( {\max}^2 \left(
   			0, v^k_i + \rho_k g_i(x^{k+1}) \right) - (v^k_i)^2 \right)
   		\\
   		&\qquad\qquad + 
   		\innerp{w^k}{h(x^{k+1})}_Y + \frac{\rho_k}{2}\nnorm{ h(x^{k+1}) }^2_Y
   		-
   		\varepsilon_k
   		\\
   		&\qquad \leq
   		f(x)
   		+
   		\frac{1}{2 \rho_k} \sum_{i=1}^m 
   			\left( {\max}^2 \left(
   			0, v^k_i + \rho_k g_i(x) \right) - (v^k_i)^2 \right)
   		\\
   		&\qquad\qquad + 
   		\innerp{w^k}{h(x)}_Y + \frac{\rho_k}{2}\nnorm{ h(x) }^2_Y
   		\end{aligned}
	\end{equation}
	for each $k\in\N$. 
	Division by $\rho_k$ yields
	\begin{equation}\label{eq:some_helpful_estimate}
		\begin{aligned}
		&\frac{f(x^{k+1})}{\rho_k} 
   		+ 
   		\frac{1}{2} \sum_{i=1}^m 
   			\left( {\max}^2 \left(
   			0, \frac{v^k_i}{\rho_k} + g_i(x^{k+1}) \right) 
   				- \left(\frac{v^k_i}{\rho_k}\right)^2 \right)
   		\\
   		&\qquad\qquad + 
   		\innerp{w^k/\rho_k}{h(x^{k+1})}_Y + \frac{1}{2}\nnorm{ h(x^{k+1}) }^2_Y
   		-
   		\frac{\varepsilon_k}{\rho_k}
   		\\
   		&\qquad \leq
   		\frac{f(x)}{\rho_k}
   		+
   		\frac{1}{2} \sum_{i=1}^m 
   			\left( {\max}^2 \left(
   			0, \frac{v^k_i}{\rho_k} + g_i(x) \right) 
   				- \left(\frac{v^k_i}{\rho_k}\right)^2 \right)
   		\\
   		&\qquad\qquad + 
   		\innerp{w^k/\rho_k}{h(x)}_Y + \frac{1}{2}\nnorm{ h(x) }^2_Y
   		\end{aligned}
   	\end{equation}
	for each $k\in\N$.
	Note that $\{f(x^{k+1})\}_{k\in K}$ is bounded from below
	due to the assumed weak sequential lower semicontinuity of $f$.
	Furthermore, we have $v^k/\rho_k\to 0$, $w^k/\rho_k\to 0$, and $\varepsilon_k/\rho_k\to 0$
	due to boundedness of $\{v^k\}$, $\{w^k\}$, and $\{\varepsilon_k\}$
	as well as $\rho_k\to\infty$. 
	Thus, taking the lower limit $k\to_K\infty$ in \eqref{eq:some_helpful_estimate} yields
	\begin{equation}\label{eq:bound_explicit_constraints}
		\frac12\norm{\max(g(\bar x),0)}^2_2+\frac12\norm{h(\bar x)}^2_Y
		\leq
		\frac12\norm{\max(g(x),0)}^2_2+\frac12\norm{h(x)}^2_Y
		<
		\infty
	\end{equation}
	where we also exploited weak sequential lower semicontinuity of $g_1,\ldots,g_m$ and
	weak-strong sequential continuity of $h$.
	Hence, $\bar x$ is a global minimizer of
	\[
		\min\limits_{x\in X} \quad \frac12\norm{\max(g(x),0)}_2^2+\frac12\norm{h(x)}_Y^2 
		\quad \text{s.t.} \quad 
		x\in \dom f\cap\dom g\cap C.
	\]
	Further, we note that for any $x\notin\dom g$, the objective value of this
	optimization problem would be $\infty$, so $\bar x$ is already a global minimizer of
	\eqref{eq:feasibility_problem}.
	
	Whenever $\dom f\cap\mathcal F$ is nonempty, 
	we may choose $x\in\dom f\cap \mathcal F$ above in order to see from
	\eqref{eq:bound_explicit_constraints} that 
	$g(\bar x)\leq 0$ and $h(\bar x)=0$ are valid, 
	and $\bar x\in\mathcal F$ follows.
	It remains to show $V_{\rho_k}(x^{k+1},v^k)\to_K 0$.
	Therefore, we exploit \cref{lem:upper_bound_AL_values} and $x\in\mathcal F$
	to obtain the estimate
	\begin{equation}\label{eq:some_really_helpful_estimate}
		\begin{aligned}
		&
		f(x^{k+1})
		- \frac{1}{2\rho_k}\nnorm{v^k}^2_2 + \innerp{w^k}{h(x^{k+1})}_Y  - \varepsilon_k
		\\		
		&\qquad 
		\leq
		f(x^{k+1}) 
   		+ 
   		\frac{1}{2 \rho_k} \sum_{i=1}^m 
   			\left( {\max}^2 \left(
   			0, v^k_i + \rho_k g_i(x^{k+1}) \right) - (v^k_i)^2 \right)
   		\\
   		&\qquad\qquad + 
   		\innerp{w^k}{h(x^{k+1})}_Y + \frac{\rho_k}{2}\nnorm{ h(x^{k+1}) }^2_Y
   		-
   		\varepsilon_k
   		\\
   		&\qquad 
   		\leq
   		f(x)
		\end{aligned}
	\end{equation}
	for each $k\in\N$ from \eqref{eq:some_interesting_estimate}.
	From $x^{k+1}\weakly_K\bar x$, we find $h(x^{k+1})\to_K h(\bar x)$ 
	by weak-strong sequential continuity of $h$.
	Thus, $\{\innerp{w^k}{h(x^{k+1})}_Y\}_{k\in K}$ remains bounded 
	as $\{w^k\}$ is bounded by construction.
	Since $\{v^k\}$ is bounded by construction as well while $\rho_k\to\infty$ holds, 
	and since $\{\varepsilon_k\}$ is assumed to be bounded, 
	$\{f(x^{k+1})\}_{k\in K}$ is bounded from above.
	Its lower boundedness has already been mentioned earlier in the proof.
	Hence, we have $f(x^{k+1})/\rho_k\to_K 0$.
	
	Dividing \eqref{eq:some_really_helpful_estimate} by $\rho_k/2$ and taking,
	in contrast to our earlier strategy, the upper limit $k\to_K\infty$
	in the second estimate gives
	\begin{align*}
		\limsup\limits_{k\to_K\infty}&
		\left(\nnorm{\max(0,v^k/\rho_k+g(x^{k+1}))}^2_2+\nnorm{h(x^{k+1})}^2_Y\right)
		\leq
		0
	\end{align*}
	as we now also know $f(x^{k+1})/\rho_k\to_K 0$.
	Due to $v^k/\rho_k\to 0$, this gives the convergence
	$\nnorm{\max(g(x^{k+1}),-v^k/\rho_k)}_2\to_K 0$,
	and $\nnorm{h(x^{k+1})}_Y\to_K 0$ follows as well.
	Since all norms in finite-dimensional spaces are equivalent,
	$V_{\rho_k}(x^{k+1},v^k)\to_K 0$ is obtained and the proof is completed.
\end{proof}

Next, we want to show that under \cref{ass:approximate_solution}, \cref{alg:ALM} can be used to
compute a global minimizer of \eqref{eq:nonsmooth_problem} provided there exists one.

\begin{theorem}\label{thm:converence_ALM_eps_optimal_solutions}
	Let $\dom f \cap \mathcal F \neq \emptyset$ hold.
	Assume that \cref{alg:ALM} produces a sequence $\{x^k\}$ such that
	\cref{ass:approximate_solution} holds for some sequence $\{\varepsilon_k\}$ 
	satisfying $\varepsilon_k\to 0$.
	Then, for each subsequence $\{x^{k+1}\}_{k\in K}$ and each point $\bar x\in X$
	satisfying $x^{k+1}\weakly_K\bar x$, we have $f(x^{k+1})\to_Kf(\bar x)$ and
	$\bar x$ is a global minimizer of \eqref{eq:nonsmooth_problem}.
\end{theorem}

\begin{proof}
	To start, note that \cref{prop:feasibility_of_accumulation_points} 
	guarantees that $\bar x$ is a feasible
	point of \eqref{eq:nonsmooth_problem}. 
	Furthermore, for each feasible point $x\in \dom f\cap\mathcal F$ 
	of \eqref{eq:nonsmooth_problem},
	\cref{ass:approximate_solution} and \cref{lem:upper_bound_AL_values} yield
	\begin{equation}\label{eq:some_estimate_global_convergence_ALM}
		\forall k\in\N\colon\quad
		L_{\rho_k}(x^{k+1},v^k,w^k)-\varepsilon_k\leq L_{\rho_k}(x,v^k,w^k)\leq f(x).
	\end{equation}
	We note that the same inequality holds trivially for all
	$x\in\mathcal F\setminus\dom f$.
	We will first prove that $\limsup_{k\to_K\infty}f(x^{k+1})\leq f(\bar x)$ is valid.
	Again, we proceed by investigating two disjoint cases.
	
	\emph{Case 1}: Suppose that $\{\rho_k\}$ remains bounded. 
		As in the proof of \cref{prop:feasibility_of_accumulation_points}, this implies
		that condition \eqref{eq:update_penalty_ALM} holds along the tail of the sequence.
		Thus, for each $i\in\{1,\ldots,m\}$, we find
	\begin{align*}
		\abs{\max\left(0,v^k_i/\rho_k+g_i(x^{k+1})\right)-v^k_i/\rho_k}
		=
		\abs{\max\left(g_i(x^{k+1}),-v^k_i/\rho_k\right)}
		\to 0
	\end{align*}
	as $k\to\infty$. By boundedness of $\{v^k_i/\rho_k\}$, $\{\max(0,v^k_i/\rho_k+g_i(x^{k+1}))\}$
	needs to be bounded as well which is why we already find
	\begin{align*}
		\abs{{\max}^2\left(0,v^k_i/\rho_k+g_i(x^{k+1})\right)-\left(v^k_i/\rho_k\right)^2}
		\to 0,
	\end{align*}
	and by boundedness of $\{\rho_k\}$, this yields
	\[
		\frac{1}{\rho_k}\left({\max}^2(0,v^k_i+\rho_k\,g_i(x^{k+1}))-(v^k_i)^2\right)
		\to 0.
	\]
	Furthermore, we find $\innerp{w^k}{h(x^{k+1})}_Y\to_K 0$ and $\tfrac{\rho_k}{2}\nnorm{h(x^{k+1})}_Y^2\to_K 0$
	from $x^{k+1}\weakly_K\bar x$, weak-strong sequential continuity of $h$, $h(\bar x)=0$, and boundedness of $\{w^k\}$.
	Plugging all this into \eqref{eq:some_estimate_global_convergence_ALM} while respecting
	the definition of the function $L_{\rho_k}$ and $\varepsilon_k\to 0$,
	we find $\limsup_{k\to_K\infty}f(x^{k+1})\leq f(x)$.
	
	\emph{Case 2}: Let $\{\rho_k\}$ be unbounded. 
	Then we already have $\rho_k\to\infty$ by construction of \cref{alg:ALM}. 
	Furthermore, \eqref{eq:some_estimate_global_convergence_ALM} implies
	validity of the estimate
	\[	
		\forall k\in\N\colon\quad
		f(x^{k+1})-\frac1{2\rho_k}\nnorm{v^k}_2^2+\innerp{w^k}{h(x^{k+1})}_Y-\varepsilon_k
		\leq
		f(x)
	\]
	by leaving out some of the nonnegative terms on the left-hand side.
	As above, we find $\innerp{w^k}{h(x^{k+1})}\to_K 0$ by boundedness of $\{w^k\}$, weak-strong sequential continuity
	of $h$, and $h(\bar x)=0$. The boundedness of $\{v^k\}$ and $\rho_k\to\infty$ yield
	$\tfrac{1}{2\rho_k}\nnorm{v^k}_2^2\to 0$ as $k\to\infty$. Thus, taking the upper limit
	in the above estimate shows $\limsup_{k\to_K\infty}f(x^{k+1})\leq f(x)$.
	
	In order to finalize the proof, 
	we observe that the weak sequential lower semicontinuity of $f$ 
	now yields the estimate
	\[
		f(\bar x)
		\leq
		\liminf\limits_{k\to_K\infty} f(x^{k+1})
		\leq
		\limsup\limits_{k\to_K\infty} f(x^{k+1})
		\leq 
		f(x).
	\]
	As this has been shown for each $x\in\dom f\cap\mathcal F$ 
	(and is trivially valid for each $x\in\mathcal F\setminus\dom f$),
	$\bar x$ is a global minimizer of \eqref{eq:nonsmooth_problem}.
	Using the above estimate with $x:=\bar x$, we additionally
	find the convergence $f(x^{k+1})\to_K f(\bar x)$.
\end{proof}

As a consequence of the previous result, we obtain the
following stronger version for convex problems with
a uniformly convex objective function.

\begin{corollary}\label{cor:converence_ALM_eps_optimal_solutions}
	Let $\mathcal F\neq\emptyset$ hold.
	Assume that \cref{alg:ALM} produces a sequence $\{x^k\}$ such that
	\cref{ass:approximate_solution} holds for some sequence $\{\varepsilon_k\}$ satisfying $\varepsilon_k\to 0$.
	Furthermore, let $f$ be continuous as well as uniformly convex, $g_1,\ldots,g_m$ be convex, $h$ be affine, and
	$C$ be convex.
	Then the entire sequence $\{x^{k}\} $ converges (strongly)
	to the uniquely determined global minimizer of 
	\eqref{eq:nonsmooth_problem}. 
\end{corollary}

\begin{proof}
Since $ f $ is uniformly convex, the (convex) optimization
problem \eqref{eq:nonsmooth_problem} has a unique solution
$ \bar x \in X $, see \cite[Theorem~2.5.1, Propositions~2.5.6, 3.5.8]{Zalinescu2002}.
As $\bar x$ is a minimizer of the underlying convex problem \eqref{eq:nonsmooth_problem}
and since $f$ is assumed to be continuous,
there exists $\bar \xi\in\partial f(\bar x)$ such that $\dual{\bar \xi}{x-\bar x}_X\geq 0$ is valid
for all $x\in\mathcal F$,
see \cite[Theorem~2.9.1]{Zalinescu2002}.
By uniform convexity of $f$, there exists a constant $ \nu > 0 $
such that 
\begin{equation}\label{eq:stronglyconvex}
	\forall x\in X\colon\quad
	f(x) \geq f(\bar x) + \langle \bar \xi, x - \bar x \rangle_X +
	\frac{\nu}{2} \| x - \bar x \|_X^2,
\end{equation}
see, e.g., the first part of the proof of \cite[Proposition~3.5.8]{Zalinescu2002}.
This implies
\begin{align*}
    f(\bar x) & + \ddual{ \bar \xi }{x^{k+1} - \bar x }_X  
	+ \frac{\nu}{2} \nnorm{ x^{k+1} - \bar x }_X^2 
	- \frac{1}{2 \rho_k} \| v^k \|_2^2 + \innerp{ w^k }{ h (x^{k+1}) }_Y \\
	& \leq f(x^{k+1}) - \frac{1}{2 \rho_k} \nnorm{ v^k }_2^2 + \innerp{ w^k }{ h (x^{k+1}) }_Y\\
	& \leq f(x^{k+1}) + \frac{1}{2 \rho_k}
	\sum_{i=1}^m \big( {\max}^2 \big\{ 0, v_i^k + \rho_k
	g_i (x^{k+1}) \} - (v_i^k)^2 \big) \\
	& \qquad + \innerp{ w^k }{ h (x^{k+1}) }_Y + \frac{\rho_k}{2} \nnorm{ h(x^{k+1}) }_Y^2 \\
    & = L_{\rho_k} (x^{k+1}, v^k, w^k) \\
    & \leq L_{\rho_k} (\bar x, v^k, w^k) + \varepsilon_k \\
    & \leq f(\bar x) + \varepsilon_k 
\end{align*}
for all $ k \in \N $, where the first inequality results from
\eqref{eq:stronglyconvex}, the second one comes from adding some
nonnegative terms, the subsequent equation is simply the
definition of the augmented Lagrangian, the penultimate
inequality takes into account \cref{ass:approximate_solution}, 
and the final estimate uses \cref{lem:upper_bound_AL_values}.

Note that, on the one hand, the term on the right-hand side is bounded. On the other
hand, since $ \{ v^k \} $ and $ \{ w^k \} $ are bounded
sequences and $ h $ is affine, the growth behavior of 
the left-hand side is dominated by the quadratic term.
Consequently, the sequence $ \{ x^k \} $ is bounded and,
therefore, has a weakly convergent subsequence in the 
reflexive space $ X $. The weak limit is necessarily
a solution of \eqref{eq:nonsmooth_problem} by 
\cref{thm:converence_ALM_eps_optimal_solutions}. Since
the entire sequence $ \{ x^k \} $ is bounded, we therefore
get $x^k\weakly\bar x$ and $ f(x^k) \to f(\bar x) $ 
from \cref{thm:converence_ALM_eps_optimal_solutions}.
 
Let us now test \eqref{eq:stronglyconvex} with $x:=x^{k+1}$. 
Then, after some rearrangements, we find
\[
	f(x^{k+1})-f(\bar x)-\dual{\bar \xi}{x^{k+1}-\bar x}_X
	\geq 
	\frac{\nu}{2}\nnorm{x^{k+1}-\bar x}^2_X
	\geq 
	0.
\]
From $x^{k+1}\weakly\bar x$ and $f(x^{k+1})\to f(\bar x)$,
the left-hand side in this estimate tends to $0$ as $k\to\infty$.
Due to $\nu>0$, this immediately gives $x^{k+1}\to\bar x$,
and the proof is complete.
\end{proof}

We end this section by discussing a suitable termination criterion for \cref{alg:ALM}
and recalling the central convergence result for convex problems stated in
\cite{Rockafellar1973} which is based on a non-safeguarded version of \cref{alg:ALM}.

\begin{remark}\label{rem:termination}
	Observe that \cref{prop:feasibility_of_accumulation_points} indicates that checking
	$V_{\rho_{k-1}}(x^k,v^{k-1})\leq\varepsilon^\textup{alm}_\textup{abs}$ for some
	$\varepsilon^\textup{alm}_\text{abs}\geq 0$ in each of the iterations $k\in\N$, $k\geq 1$, 
	is a reasonable termination criterion for \cref{alg:ALM}
	provided $\dom f\cap\mathcal F\neq\emptyset$.
	On the one hand, if $V_{\rho_{k-1}}(x^k,v^{k-1})$ is small, then the underlying point
	$x^k$ is close to be feasible, and the associated Lagrange multiplier estimate $v^{k-1}$
	is close to satisfy the associated complementarity-slackness condition w.r.t.\ the
	inequality constraints.
	On the other hand, along weakly convergent subsequences of the iterates 
	produced by \cref{alg:ALM}, $V_{\rho_{k-1}}(x^k,v^{k-1})$ indeed becomes arbitrarily small
	under the assumptions of \cref{prop:feasibility_of_accumulation_points}.
	Furthermore, under the assumptions of \cref{thm:converence_ALM_eps_optimal_solutions}, 
	weak accumulation points are already global minimizers of \eqref{eq:nonsmooth_problem}.
\end{remark}

\begin{remark}\label{rem:rockafellar_convex_problems}
	In \cite{Rockafellar1973}, the author considers \eqref{eq:nonsmooth_problem} with
	convex, continuously differentiable functions $f,g_1,\ldots,g_m$, a closed, convex
	set $C$, and in the absence of (infinite-dimensional) equality constraints.
	For the numerical solution of this problem, a simplified version of \cref{alg:ALM}
	is suggested where safeguarding of the Lagrange multipliers is omitted and the
	penalty parameter stays constant.
	More precisely, this means that 
	\cref{item:choice_of_approximate_multipliers,item:update_penalty_ALM} are removed from
	the method, in \cref{item:solve:subproblem_ALM}, the subproblem
	\begin{equation}\label{eq:simplified_subproblem}
		\min_{x\in X} \quad L_{\rho_0} (x, \lambda^k) \quad \text{s.t.} \quad  x \in C.
    \end{equation}
    is solved up to $\varepsilon_k$-minimality, see \cref{ass:approximate_solution}, 
    and the Lagrange multiplier update
    in \cref{item:adjust_multiplier_ALM} is replaced by 
    \[
    	\lambda^{k+1}:=\max\bigl(0,\lambda^k+\rho_0g(x^{k+1})\bigr).
	\]
	In \cite[Theorem~2.1]{Rockafellar1973}, it is shown that whenever \eqref{eq:nonsmooth_problem}
	possesses a minimizer which is stationary (in the sense that Lagrange multipliers exist)
	while $\sum_{k=0}^\infty\sqrt{\varepsilon_k}<\infty$ holds, then each accumulation point
	$(\bar x,\bar\lambda)\in X\times\R^m_+$ of the primal-dual sequence 
	$\{(x^k,\lambda^k)\}$ produced by the adjusted method is a Karush--Kuhn--Tucker pair 
	of \eqref{eq:nonsmooth_problem}, i.e., $\bar x$ is a feasible point with associated
	Lagrange multipliers $\bar\lambda$.
	Particularly, $\bar x$ is a minimizer of \eqref{eq:nonsmooth_problem} in this case.
\end{remark}

\section{Numerical solution of the Poisson denoising model}\label{sec:experiments}
In this section, we numerically test the variational Poisson denoising model. 
First, in \cref{sec:reformulation}, we introduce an affine reformulation of the constraints 
in \eqref{eq:VPD}. 
General comments about the implementation of the experiments
are presented in \cref{sec:denoising_implementation}.
Particularly, we introduce three augmented Lagrangian schemes which are used
to tackle \eqref{eq:VPD} computationally, one of them focusing on the original model
\eqref{eq:VPD} while the other two address the aforementioned reformulation.
In \cref{sec:denoising_SGD}, we comment on a stochastic gradient descent method 
which is used to solve the arising subproblems within the augmented Lagrangian schemes. 
Our way of documenting the obtained results is briefly outlined in \cref{sec:documentation}.
Finally, the outcome of our experiments is presented in \cref{sec:denoising_numerical_results}.

\subsection{Reformulation of the constraints}\label{sec:reformulation}

Let us point out that the constraints in \eqref{eq:VPD} are heavily nonlinear and,
due to the fact that the Kullback--Leibler divergence is extended real-valued,
nonsmooth. In this subsection, we show that these constraints can be reformulated as
affine ones, making the associated reformulation of \eqref{eq:VPD} seemingly
easier to tackle from an algorithmic point of view.

For given $a\geq 0$ and $c>0$, let us consider the nonlinear constraint
\begin{equation}\label{eq:nonlinear_constraint}
	\eta(a,b)\leq c
\end{equation}
on the real variable $b$.
Above, $\eta$ is the Kullback--Leibler divergence defined in \eqref{eq:Kullback_Leibler_Divergence}.
Whenever $a=0$, we can immediately rewrite \eqref{eq:nonlinear_constraint} as $0\leq b\leq c$.

Thus, let us assume that $a>0$.
Then \eqref{eq:nonlinear_constraint} requires $b>0$.
Let us study the properties of the smooth function $(0,\infty)\ni b\mapsto\eta(a,b)\in\R$.
One can easily check that it is strictly convex with a uniquely determined minimizer at
$\check b:=a$ with function value $0$.
Furthermore, we have
\[
	\lim\limits_{b\downarrow 0}\eta(a,b)=\infty,
	\qquad
	\lim\limits_{b\to\infty}\eta(a,b)=\infty.
\]
Hence, there exist uniquely determined points $\underline b\in(0,\check b)$ and
$\overline b\in(\check b,\infty)$ which solve the nonlinear equation
\begin{equation}\label{eq:nonlinear_equation}
	\eta(a,b)-c=0,
\end{equation}
and, thus,
\eqref{eq:nonlinear_constraint} is equivalent to $\underline b \leq b\leq \overline b$.
Let us now describe how the bounds $\underline b$ and $\overline b$ can be accessed computationally.
Division by $a$ in \eqref{eq:nonlinear_equation} and some rearrangements
yield
\[	
	\frac{b}{a} + \ln\left(\frac{a}{b}\right) = \frac{c}{a}+1.
\]
We exponentiate this equation and take the reciprocal on both sides in order to find
\[
	\left(-\frac{b}{a}\right)\exp\left(-\frac{b}{a}\right)
	=
	\frac{-1}{\exp\left(\frac{c}{a}+1\right)}.
\]
Due to
\[
	\frac{-1}{\exp(1)}
	<
	\frac{-1}{\exp\left(\frac{c}{a}+1\right)}
	<
	0,
\]
the latter equation indeed possesses two solutions which can be expressed
in terms of \emph{Lambert's $W$ function} $W_\kappa$, 
which is the multivalued inverse of $t\mapsto t\exp(t)$,
see e.g.\ \cite{CorlessGonnetHareJeffreyKnuth1996}
for properties, applications, and the numerical evaluation of this function,
as 
\begin{align*}
	\underline{b}:=\underline{b}(a,c):=-a W_0\left(\frac{-1}{\exp\left(\frac{c}{a}+1\right)}\right),
	\quad
	\overline{b}:=\overline{b}(a,c):=-a W_{-1}\left(\frac{-1}{\exp\left(\frac{c}{a}+1\right)}\right).
\end{align*}

Taking all our above findings together, we can equivalently state \eqref{eq:VPD} as
\begin{equation}\label{eq:VPDaff}\tag{VPD$_\textup{aff}$}
	\begin{aligned}
		&\min\limits_{u\in L^2(\Omega)}&	&f(u)&	&&
		\\
		&\text{ s.t. }&	&0\leq u_B\leq r(|B|)&	&\forall B\in\mathcal B_Z^0&
		\\
		&&				&\underline{b}(Z_B,r(|B|))\leq u_B\leq \overline{b}(Z_B,r(|B|))&	
			&\forall B\in\mathcal B_Z^+,&
	\end{aligned}
\end{equation}
where we made use of the sets
\[
	\mathcal B_Z^0:=\{B\in\mathcal B\,|\,Z_B=0\},
	\qquad
	\mathcal B_Z^+:=\{B\in\mathcal B\,|\,Z_B>0\}.
\]
Observe that \eqref{eq:VPDaff} possesses twice as many constraints as \eqref{eq:VPD} does.

Let us note that \eqref{eq:VPDaff} is \emph{not} a box-constrained problem since $u_B$ already
involves an averaging operation of $u$ on box $B\in\mathcal B$.
However, using the weighted characteristic function $\overline\chi_B\in L^2(\Omega)$
given by
\[
	\forall \omega\in\Omega\colon\quad
	\overline\chi_B(\omega)
	:=
	\begin{cases}
		|B|^{-1}	&	\omega\in B,\\
		0			&	\omega\notin B
	\end{cases}
\]
for each subbox $B\in\mathcal B$,
the constraints in \eqref{eq:VPDaff} can be rewritten in the form
\[
	\begin{aligned}
		&0\leq \dual{\overline\chi_B}{u}_{L^2(\Omega)}\leq r(|B|)&	&\forall B\in\mathcal B_Z^0&
		\\
		&\underline{b}(Z_B,r(|B|))\leq \dual{\overline\chi_B}{u}_{L^2(\Omega)}\leq \overline{b}(Z_B,r(|B|))&	
			&\forall B\in\mathcal B_Z^+.&
	\end{aligned}
\]
Hence, whenever the objective function $f$ is smooth enough, 
local minimizers of \eqref{eq:VPDaff} are stationary point of that problem
due to \cite[Proposition~2.42]{BonnansShapiro2000}, i.e., Lagrange multipliers exist
in that case.

\subsection{Implementation}\label{sec:denoising_implementation}

For the numerical realization, we discretize the image using $256^2$ equally sized pixels. 
We use a zero padding of five equally sized pixels at the borders of the image to prevent artifacts in the reconstruction 
due to the periodic extension in the objective value associated with \eqref{eq:sobolev_penalty}. 
Thus, setting $n:=266$, the image to be denoised has actually a size of $n^2$ pixels. 
A wider padding would give a higher guarantee of avoiding these artifacts at the cost of higher computational complexity. 
Empirically, our tests show that the padding size of five pixels seems to be sufficient.  
The (padded) image $u$ 
is therefore approximated by an $n\times n$ matrix of pixels with pixel-size $s := 1/n^2 = 266^{-2}$.
With this resolution, the family of regions $\mathcal{B} \subset 2^{\Omega}$ is chosen as all 
subsquares of the image with side length (scale) between $1$ and $64$ pixels. 
The size $|B|$ of a region $B\in\mathcal{B}$ is numerically computed as $|B| := s\, \#B$.
The number of constraints in \eqref{eq:VPD}, i.e., the number of subboxes in $\mathcal B$, is then calculated by
\begin{align*}
	\sum_{i = n-63}^{n} i^2 = 3541216,
\end{align*}
and approximately $50$ times higher than the total number of pixels.
The usage of squares is clearly subjective, but without a priori information
about the image content, it is not clear which other shapes might be preferable.
We emphasize that, however, our methodology is also applicable with other shapes
such as rectangles, circles, or ellipses.

As there are way more subsquares with small side length, a penalty term
\[
\operatorname{pen}(|B|) := \sqrt{2 \left( \log\left(n^2/\left| B \right|\right)+1 \right)},
\]
which only depends on the size of the subsquares, is introduced.
This is necessary to avoid the small subsquares to dominate the statistical behavior of the overall test statistic
\[
T_n(Z,u,\mathcal{B}) := \max_{B \in \mathcal{B}}[T_B(Z,u) - \operatorname{pen}(|B|)],
\]
see \cite{KoenigMunkWerner2020},
where $T_B$ is the LRT statistic from \eqref{eq:LRT_statistic}.
We approximate the $(1-\alpha)$-quantile $q_{1-\alpha}$ of $T_n$ by the (empirically sampled) $(1-\alpha)$-quantile $\tilde{q}_{1-\alpha}$ of
\[
M_n(\mathcal{B}) := \max_{B \in \mathcal{B}} \left[ |B|^{-1/2} \left| \sum\nolimits_{i \in B} X_i \right| - \operatorname{pen}(|B|) \right]
\]
with i.i.d.\ standard normal random variables $X_i$.
If the smallest scale in $\mathcal B$ was at least of size $\log(n)$, then this approximation was shown to be valid in \cite{KoenigMunkWerner2020}. 
However, the chosen penalization $\operatorname{pen}$ effectively overdamps the small scales, 
see \cite{SharpnackAriasCastro2016}, 
which makes this approximation reasonable over all scales considered here. 
Altogether, this leads to the right-hand side
\[
r(|B|) := \frac{\left(\tilde{q}_{1-\alpha}+\operatorname{pen}(|B|)\right)^2}{2|B|}
\]
in \eqref{eq:VPD}. 
In the numerical experiments, the $0.1$-quantile $\tilde{q}_{0.1} := 1.63$ is used, 
because for bigger values of $\alpha$, the local hypothesis tests are not restrictive 
enough such that the obtained reconstruction is oversmoothed. 
For the same reason, $\theta:=0.01$ is chosen relatively small in the Sobolev-type penalty \eqref{eq:sobolev_penalty},
which is used as the objective function subsequently.
Observe that the discretized problem associated with \eqref{eq:VPDaff} is a convex quadratic optimization problem
in this situation.

We solved the associated problem \eqref{eq:VPD} with the aid of the following procedures:
\begin{itemize}[leftmargin=6em]
\item[\ALM:] 
	the safeguarded augmented Lagrangian method from \cref{alg:ALM} applied to \eqref{eq:VPD},	
\item[\ALMR:] 
	the safeguarded augmented Lagrangian method from \cref{alg:ALM} applied to the
	reformulated problem \eqref{eq:VPDaff}, and
\item[\SALMR:] 
	the simplified augmented Lagrangian method without safeguarding and
	with constant penalty parameter applied to \eqref{eq:VPDaff},
	see \cref{rem:rockafellar_convex_problems}.
\end{itemize}
For all three methods, the discretized noisy observation is taken as the primal starting point, i.e., $x^0:=Z$,
while $\lambda^0:=0$ is chosen for the initial Lagrange multiplier.
Note that this choice for $x^0$ is only possible in the discretized setting,
as $Z$ is likely to lack $L^2$-regularity in the infinite-dimensional framework. 
In the case of continuous computations, one could e.g.\ use a kernel density
estimator to obtain some point $x^0 \in L^2 (\Omega)$ from $Z$.
Within our implementation of \ALM\ and \ALMR,
we choose $v^k$ as the componentwise projection
of the Lagrange multiplier $\lambda^k$ onto the interval $[0,10^8]$. 
Furthermore, we use the parameters 
$\rho_0 := 4$, $\tau := 0.9$, and $\gamma := 4$ for these methods.
In \SALMR, we make use of the (constant) penalty parameter $\rho_0:=4\cdot 10^5$.
Finally, for the abort of all three methods,
we exploit the termination criterion from \cref{rem:termination} 
with $\varepsilon_{\textup{abs}}^{\textup{alm}} := 10^{-2}$. 
Let us note that the aforementioned choice for $\rho_0$ in \SALMR\ 
is such that this termination criterion is almost satisfied already 
in the first iteration of the method.
In some preliminary experiments, it turned out that \SALMR\ is not capable
to enhance the feasibility properties along the iterates significantly
if $\rho_0$ is chosen much smaller than $4\cdot 10^{5}$.
Finally, all three methods are, at the latest, aborted after a total number of $30$ outer
augmented Lagrangian iterations.

\subsection{Stochastic gradient descent as a subproblem solver}\label{sec:denoising_SGD}

Solving the unconstrained associated subproblems \eqref{eq:ALM_subproblem} as well as \eqref{eq:simplified_subproblem}
within the algorithmic frameworks \ALM\ as well as \ALMR\ and \SALMR, respectively,
is computationally expensive,
especially as the problems \eqref{eq:VPD} and \eqref{eq:VPDaff} possess many constraints. 
This obstacle is tackled by using the first-order gradient descent method NADAM from \cite{Dozat2016} 
which outperformed other gradient descent methods 
in the exemplary setting we are considering here. 
It is also utilized that the constraints are redundant to a certain degree and, thus, 
a stochastic version of the NADAM method can be used.
Undoubtedly, this also speeds up the evaluation of the (stochastic) gradient of the
augmented Lagrangian function and cheapens each iteration of the subproblem solver. 

Let us comment on the implementation of NADAM in the context of \ALM.
For the fixed penalty parameter $\rho_k>0$ and the Lagrange multiplier estimate $v^k:=\{v^k_B\}_{B\in\mathcal B}$, 
the augmented Lagrangian subproblem \eqref{eq:ALM_subproblem} takes the particular form
\begin{equation*}
	\min\limits_{u\in\R^{n\times n}}
	\quad
	f(u)+\frac1{2\rho_k}\sum\limits_{B\in\mathcal B}
		\left(
			\max{}^2\bigl(0,v^k_B+\rho_k(\eta(Z_B,u_B)-r(|B|))\bigr)-(v_B^k)^2
		\right)
\end{equation*}
in the present situation. Here and in what follows, we approximate the continuous mean 
$|B|^{-1} \int_B u(\omega)\, \mathrm d \omega$ by $u_B := s|B|^{-1} \sum_{i \in B} u_i = (\# B)^{-1} \sum_{i \in B} u_i$, 
which corresponds to the discrete mean. 
For the NADAM method, one needs to calculate the gradient of the augmented Lagrangian function. 
Therefore, the partial derivative of $u\mapsto \eta(Z_B,u_B)$ w.r.t.\ pixel $u_i$ (where it exists) is given by
$(\# B)^{-1}(1-Z_B/u_B)$ if $i\in B$ and, otherwise, $0$.
Thus, the partial derivative of the associated augmented Lagrangian function w.r.t.\ the pixel $u_i$ (where it exists) equals
\[
	f'_{u_i}(u)
	+ 
	\sum_{B \in \mathcal{B}(i) } 
		\frac{1}{\# B}\max\left(0, v^k_B + \rho_k(\eta(Z_B, u_B) - r(|B|))\right)\left(1- \frac{Z_B}{u_B} \right),
\]
where
\begin{equation*}
	\mathcal{B}(i) := \{B \in \mathcal{B} \,|\, i \in B\}.
\end{equation*}
The above formula is valid whenever $u_B > 0$ for all $B \in \mathcal B(i)$. To account for the non-differentiability on the boundary, we set
\[
	(L_{\rho_k})'_{u_i}(u,v^k)
	:= 
	f'_{u_i}(u)
	+ 
	\sum_{B \in \mathcal{B}(i)}  b_{\rho_k}(Z, u, B, v^k)
\]
with 
\begin{align*}
	&b_{\rho_k}(Z, u, B, v^k)
	\\
	&:= 
	\begin{cases} 
	\frac{1}{\# B}\max\left(0, v^k_B + \rho_k(\eta(Z_B, u_B) - r(|B|))\right)\left(1- \frac{Z_B}{u_B} \right) 	& \text{if } u_B > 0, \\
	C								& \text{if } u_B = 0 \text{ and } Z_B > 0\\
	0								& \text{if } u_B = Z_B = 0.
	\end{cases}
\end{align*}
In the case $u_B = Z_B = 0$, the constraint is satisfied and, thus, we can set the corresponding gradient to $0$.
The rationale behind the definition for $u_B = 0$ and $Z_B > 0$ is to enforce a step into positive direction.
Numerically, we use $C := -10 < 0$.
By definition it always holds $Z_B \geq 0$ and, thus, we do not need to cope with the case $Z_B < 0$. 
As the NADAM method may produce iterates having pixels with negative value, 
we set all pixels to zero which have negative value after one NADAM iteration. 
This way, we also ensure $u_B \geq 0$.  

Instead of calculating the summand for every $B \in \mathcal{B}$, 
we choose a random family $\mathcal{B}_r \subset \mathcal{B}$ and approximate the gradient by
\[
	(L_{\rho_k})'_{u_i}(u,v^k)
	\approx 
	f'_{u_i}(u) 
	+ 
	\sum_{B \in \mathcal{B}_r \cap \mathcal{B}(i)} b_{\rho_k}(Z, u, B, v^k).
\]
As it is possible to efficiently calculate all summands with same scale $|B|$ with the help of the discrete Fourier transform, 
we pick only the scales at random and include all sets $B$ of those scales in $\mathcal{B}_r$.
In practice, it was first tried to use a fixed number of $10$ scales. 
This yielded fast convergence in the beginning, but convergence slowed down during the runs due to missing accuracy in solving the subproblems. 
Thus, we decided to increase the number of scales picked during the algorithm although this worsens the running time of a single augmented Lagrangian step. 
More precisely, in our experiments, we now increase the amount of scales by one after every augmented Lagrangian step. 
For simplicity, a fixed number of $300$ iterations of the NADAM method is chosen, 
and the stepsize is picked constant as $\max(0.005, 0.8^k)$ in the $k$-th iteration of \cref{alg:ALM} 
to solve the augmented Lagrangian subproblem \eqref{eq:ALM_subproblem}.

Within the algorithmic framework of \ALMR, 
for given penalty parameter $\rho_k>0$ and Lagrange multiplier estimates $v^k:=\{\underline{v}^k_B\}_{B\in\mathcal B}\cup\{\overline{v}^k_B\}_{B\in\mathcal B}$,
the associated augmented Lagrangian subproblem \eqref{eq:ALM_subproblem} is given by
\begin{equation}\label{eq:subproblem_ALMR}
	\begin{aligned}
	\min\limits_{u\in\R^{n\times n}}
	\quad
	f(u)
	&+
	\frac1{2\rho_k}\sum\limits_{B\in\mathcal B_Z^0}
		\Bigl(
			\max{}^2\bigl(0,\overline{v}^k_B+\rho_k(u_B-r(|B|))\bigr)-(\overline{v}_B^k)^2
		\\
		&\qquad\qquad\qquad
			+
			\max{}^2\bigl(0,\underline{v}^k_B-\rho_k u_B\bigr)-(\underline{v}_B^k)^2
		\Bigr)
	\\
	&+
	\frac1{2\rho_k}\sum\limits_{B\in\mathcal B_Z^+}
		\Bigl(
			\max{}^2\bigl(0,\overline{v}^k_B+\rho_k(u_B-\overline{b}(Z_B,r(|B|)))\bigr)-(\overline{v}_B^k)^2
		\\
		&\qquad\qquad\qquad
			+
			\max{}^2\bigl(0,\underline{v}^k_B+\rho_k(\underline{b}(Z_B,r(|B|))-u_B)\bigr)-(\underline{v}_B^k)^2
		\Bigr).
	\end{aligned}
\end{equation}
In the discretized setting, the partial derivative of the augmented Lagrangian function w.r.t.\ the pixel $u_i$ equals
\begin{align*}
	f'_{u_i}(u)
	&+
	\sum\limits_{B\in\mathcal B_Z^0\cap\mathcal B(i)}
		\frac{1}{\# B}\Bigl(
			\max\bigl(0,\overline{v}^k_B+\rho_k(u_B-r(|B|))\bigr)-\max\bigl(0,\underline{v}^k_B-\rho_k u_B\bigr)
		\Bigr)
	\\
	&+
	\sum\limits_{B\in\mathcal B_Z^+\cap\mathcal B(i)}
		\frac{1}{\# B}\Bigl(
			\max\bigl(0,\overline{v}^k_B+\rho_k(u_B-\overline{b}(Z_B,r(|B|)))\bigr)
			\\
		&\qquad\qquad\qquad\qquad
			-\max\bigl(0,\underline{v}^k_B+\rho_k(\underline{b}(Z_B,r(|B|))-u_B)\bigr)
		\Bigr). 
\end{align*}
Again, we do not exploit the full gradient of the augmented Lagrangian function in the NADAM framework
but rely on an approximation where the appearing sums are restricted to a random family
$\mathcal B_r\subset\mathcal B$ which is chosen in the same way as outlined above.
Similarly, the maximum number of iterations and the stepsize for NADAM are chosen 
as described in the setting of \ALM.

Finally, in the setting of \SALMR, 
the appearing augmented Lagrangian subproblem \eqref{eq:simplified_subproblem}
equals \eqref{eq:subproblem_ALMR} with $\rho_k:=\rho_0$ and $v^k:=\lambda^k$ for the current
Lagrange multiplier $\lambda^k:=\{\underline{\lambda}^k_B\}_{B\in\mathcal B}\cup\{\overline{\lambda}^k_B\}_{B\in\mathcal B}$,
and this subproblem is treated in the same way 
as \eqref{eq:subproblem_ALMR} in the \ALMR\ framework. 
The only exception is that we choose the number of NADAM iterations by $10^4$ 
in the first outer iteration of the augmented Lagrangian method. 
This is necessary to get a good approximate solution of the subproblem. 
In the subsequent outer iterations, we again use merely $300$ NADAM  iterations
as we are already closer to a good approximate solution.  

\subsection{Documentation of results}\label{sec:documentation}

In order to compare the algorithmic frameworks in $\mathcal A:=\{\ALM,\ALMR,\SALMR\}$ in a reasonable way,
we run all these methods on the same set $\mathcal P$ of benchmark instances for \eqref{eq:VPD} and illustrate
the outcome with the aid of so-called performance profiles, see \cite{DolanMore2002},
based on different performance indicators.

Let us briefly comment on the concept of performance profiles.
Therefore, we denote the output of algorithm $\mathrm a \in\mathcal A$
for problem $\mathrm p\in\mathcal P$ by $u_{\mathrm p}^{\mathrm a}$.
Let $q$ be a given performance measure taking only positive values
such that a smaller value of $q$ indicates a better performance.
The associated performance ratio is defined via
\[
	\forall \mathrm p\in\mathcal P,\,\forall \mathrm a\in\mathcal A\colon\quad
	r_{\mathrm{p},\mathrm{a}}
	:=
	\frac{q(u_{\mathrm p}^{\mathrm a})}{\min\{q(u_{\mathrm p}^{\mathrm a'})\,|\,\mathrm a'\in\mathcal A\}}.
\]
In the associated performance profile, for each of the algorithms $\mathrm a\in\mathcal A$,
we plot the illustrative part of the nondecreasing curve $\zeta_{\mathrm a}\colon[1,\infty)\to[0,1]$ given by
\[
	\forall t\in[1,\infty)\colon\quad
	\zeta_{\mathrm a}(t)
	:=
	\frac{\#\{\mathrm p\in\mathcal P\,|\,r_{\mathrm p,\mathrm a}\leq t\}}{\#\mathcal P}.
\]
Exemplary, for $\mathrm a\in\mathcal A$, $\zeta_{\mathrm a}(1)$ denotes the portion of benchmark problems
in $\mathcal P$ on which algorithm $\mathrm a$ performed best.
For $t>1$, $\zeta_{\mathrm a}(t)$ represents the portion of benchmark problems $\mathrm p\in\mathcal P$
on which algorithm $\mathrm a$ produces an output whose value w.r.t.\ the performance measure $q(u_{\mathrm p}^{\mathrm a})$ is at
most the $t$-fold of the best value achieved by any algorithm in $\mathcal A$ for the particular instance $\mathrm p$.

To document the results of our experiments, we make use of the subsequently elucidated performance measures.
\begin{itemize}
	\item \underline{\emph{Objective function value of output}}\\
		We note that the Sobolev-type penalty $f$ from \eqref{eq:sobolev_penalty} only takes nonnegative values.
		In all practically relevant scenarios, all feasible points of \eqref{eq:VPD} will possess positive
		objective value for this choice of the objective function.
	\item \underline{\emph{Number of outer augmented Lagrangian iterations until termination}}
	\item \underline{\emph{Value of error measure when the method is aborted}}\\
		For $\mathrm a\in\{\ALM,\ALMR\}$, the performance measure is given by $V_{\rho_{k-1}}(u_{\mathrm p}^{\mathrm a},v^{k-1})$,
		where $k\in\N$ is the final value of the iteration counter when \cref{alg:ALM} is aborted.
		For $\mathrm a:=\SALMR$, we take $V_{\rho_0}(u_{\mathrm p}^{\mathrm a},\lambda^{k-1})$ as the performance measure,
		where the meaning of $k$ is the same as mentioned above.
	\item \underline{\emph{Value of the penalty parameter when the method is aborted}}\\
		Here, we only compare \ALM\ and \ALMR\ since the penalty parameter is not enlarged in \SALMR.
	\item \underline{\emph{Percentage of constraints violated by the output}}\\
		Let us note that a constraint in \eqref{eq:VPD} is violated if and only if precisely one of the
		associated two constraints in \eqref{eq:VPDaff} is violated.
		Particularly, an arbitrarily chosen point satisfies at least 50\% of the constraints in \eqref{eq:VPDaff}.
		We, thus, measure the percentage of violated constraints in terms of the constraints in \eqref{eq:VPD}
		for all three methods in order to state a fair comparison.
		Due to the huge number of constraints, it is unlikely that any of the methods 
		actually computes a feasible point of \eqref{eq:VPD} or \eqref{eq:VPDaff}.
		Nevertheless, we added the positive offset $\delta:=1$ to each of the percentages 
		to obtain meaningful performance profiles.
	\item \underline{\emph{Maximum relative (affine) constraint violation of the output}}\\
		For a comparatively fair comparison, this quantity is first measured in terms of the constraints of \eqref{eq:VPD} only, i.e.,
		for any $\mathrm p\in\mathcal P$ and $\mathrm a\in\mathcal A$, the quantity
		\[
			\max\limits_{B\in\mathcal B}\frac{\max(0,\eta(Z_B,(u_{\mathrm p}^{\mathrm a})_B)-r(|B|))}{r(|B|)}
		\]
		is used as a performance measure. 
		For the outputs of \ALM, this quantity has almost always been finite.
		However, as \ALMR\ and \SALMR\ aim to compute feasible points of \eqref{eq:VPDaff},
		it may happen that, for some boxes $B\in\mathcal B$, the appearing fraction takes value $\infty$, 
		see \eqref{eq:Kullback_Leibler_Divergence} again for the precise definition of the
		Kullback--Leibler divergence $\eta$.
		In this case, we set the maximum relative constraint violation to be the largest finite value of
		\[
			\frac{\max(0,\eta(Z_B,(u_{\mathrm p'}^{\mathrm a'})_B)-r(|B|))}{r(|B|)}
		\]
		which is observed for all $\mathrm p'\in\mathcal P$ and $\mathrm a'\in\mathcal A$
		in order to get a meaningful performance profile.
		Second, in similar fashion, the relative constraint violation is measured in terms of the
		constraints in \eqref{eq:VPDaff} for all three methods.
		Note that, in this case, no issues regarding infinite values have to be faced.
		Again, a constant offset of $\delta:=1$ is added to the relative (affine) constraint violation
		to make the associated performance profile compelling.
	\item \underline{\emph{Average relative (affine) constraint violation of the output}}\\
		This performance measure is first defined by
		\[
			\sum_{B\in\mathcal B}\frac{\max(0,\eta(Z_B,(u_{\mathrm p}^{\mathrm a})_B)-r(|B|))}{\#B\,r(|B|)}
		\]
		for $\mathrm p\in\mathcal P$ and $\mathrm a\in\mathcal A$.
		Again, we note that the appearing fraction may take value $\infty$ for $\mathrm a\in\{\ALMR,\SALMR\}$.
		In this case, it is replaced by the largest finite value of this fraction achieved over $\mathcal P$ and $\mathcal A$.
		Second, a similar quantity based on the affine constraints in \eqref{eq:VPDaff} is defined for a comparison
		of all three methods.
		In both situations, we use the constant offset $\delta:=1$ to obtain meaningful performance profiles.
\end{itemize}

Let us mention that we do not rely on computation time as a performance measure for the following reasons.
Whenever $k\in\N$ denotes the final value of the outer iteration counter when the augmented Lagrangian
method of interest is aborted,
then, by construction, \ALM\ and \ALMR\ have run precisely $300k$ iterations of the subproblem solver NADAM
while \SALMR\ has run precisely $10^4+300(k-1)\approx 300(k+32)$ iterations of it.
Let us also emphasize that, with increasing outer iteration counter, the cost for one iteration of NADAM
increases as well since the number of subboxes, taken into account for the computation of the stochastic gradient
of the augmented Lagrangian function, is enlarged in each outer iteration,
see \cref{sec:denoising_SGD}.
Thus, $k$ is the decisive quantity for the running time of all three methods, 
and \SALMR\ is, generally, slower than \ALM\ and \ALMR\ since the maximum number of outer iterations is set,
for all three methods, to $30$, see \cref{sec:denoising_implementation}.

\subsection{Numerical results}\label{sec:denoising_numerical_results}

\begin{figure}[!htb]
	\centering
	{%
		\includegraphics[width=0.325\linewidth]{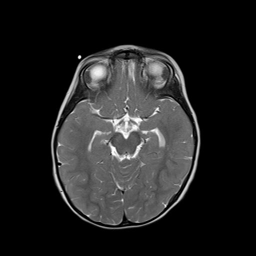}
		\includegraphics[width=0.325\linewidth]{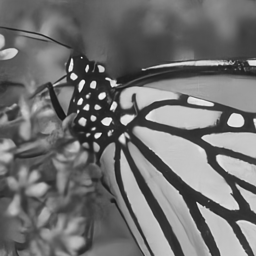}
		\includegraphics[width=0.325\linewidth]{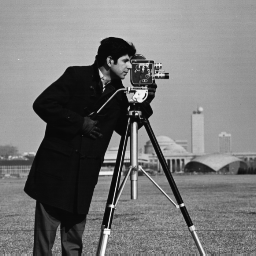}
	}%
	\\
	{%
		\includegraphics[width=0.325\linewidth]{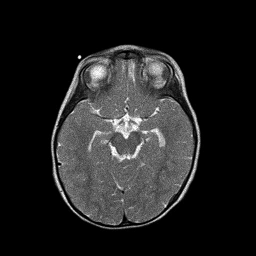}
		\includegraphics[width=0.325\linewidth]{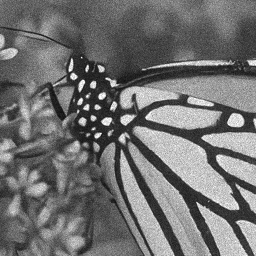}
		\includegraphics[width=0.325\linewidth]{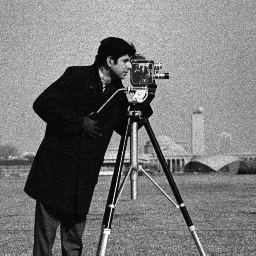}
	}%
	\\
	\caption{%
		The three test images (top row) 
		and noisy versions of the test images (bottom row).
		}%
	\label{fig:test_images}
\end{figure}

In this section, we comment on the numerical behavior of the 
augmented Lagrangian schemes \ALM, \ALMR, and \SALMR\
for the denoising of the three standard test images 
``Brain'', ``Butterfly'', and ``Cameraman'',
where the hyper-parameters are chosen as described in the previous sections. 
The benchmark problems are created by applying Poisson noise 
to each of the three ground truth pictures $\hat u$. 
Thus, every pixel $Z_i$ of the (discrete) noisy observation $Z$ is randomly drawn according to 
\begin{align*}
	\forall i \in \{1,\dots,n\}^2\colon\quad
	Z_{i}\sim \Poi(\hat u_{i}).
\end{align*}
To obtain a reasonable benchmark collection of noisy test images,
we created 10 noisy versions of each of the three test images
mentioned above.
The three test images of interest as well as some exemplary noisy versions
of them can be found in \cref{fig:test_images}.

\begin{figure}[!htb]
	\centering
	{%
		\includegraphics[width=0.325\linewidth]{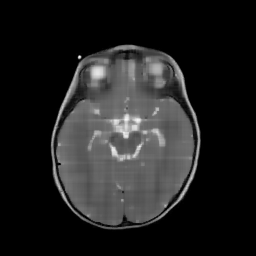}
		\includegraphics[width=0.325\linewidth]{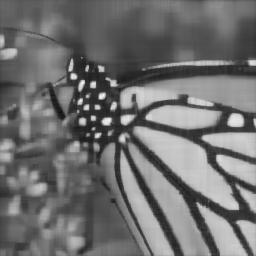}
		\includegraphics[width=0.325\linewidth]{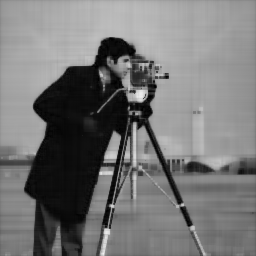}
	}%
	\\
	\caption{%
		Denoised versions of the test images obtained via \ALM.
		}%
	\label{fig:denoised_images}
\end{figure}

We pick $r$ such that \eqref{eq:quantile} and, consequently, \eqref{eq:guarantee} hold true with $\alpha = 0.1$.
Typical examples of outputs associated with \ALM\ can be found in \cref{fig:denoised_images}.
The reconstructions show that the method yields reasonable results as convergence to a meaningful solution is observed. 
As mentioned before, using $\alpha = 0.1$ ensures the statistical guarantee \eqref{eq:guarantee}, 
and, thus, leads by construction to a method tending to oversmoothing. 
This is clearly visible in \cref{fig:denoised_images}, 
but must be seen as a feature of the variational Poisson denoising method under consideration. 
Let us note that the outputs of \ALMR\ and \SALMR\ are of comparable quality 
so we do not present the associated denoised images here. 

\begin{figure}
\centering
\includegraphics[width=0.35\textwidth]{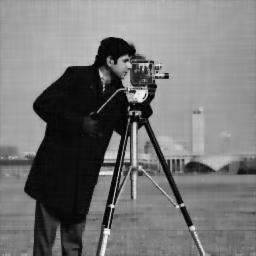}
\caption{%
	Reconstruction of the ``cameraman'' image with decreased $r$ 
	to reduce oversmoothing obtained via \ALM.
}%
\label{fig:cameraman}
\end{figure}

To cope with this, one could apply the variational Poisson denoising approach 
with a smaller function $r$ (shifted by a constant) 
to prevent oversmoothing.
Exemplary, this can be seen in the denoised version of the ``cameraman'' image
in \cref{fig:cameraman} which, again, has been
obtained with the aid of \ALM\ but with decreased $r$.
This modification implies that the statistical guarantee \eqref{eq:guarantee} is lost, 
but therefore the constraints are more restrictive and prevent oversmoothing.
Note that a similar observation was made in \cite{FrickMarnitzMunk2013}. 
It can be seen from \cref{fig:cameraman} that the corresponding reconstruction is less smooth, 
but seems visually superior over the one in \cref{fig:denoised_images}. 
For our numerical comparison of the three methods
\ALM, \ALMR, and \SALMR, however, we stick to our original choice of $r$ and
$\alpha$ in order to preserve the idea behind variational Poisson denoising.

Let us now start to document the comparison of the three augmented Lagrangian methods
\ALM, \ALMR, and \SALMR\ based on the performance measures mentioned in \cref{sec:documentation}.
Let us first inspect the performance profiles in \cref{fig:basic_PPs}
which document the obtained objective function values,
the total number of outer augmented Lagrangian iterations, as well as the value of the error measure and
the size of the penalty parameter when the respective method is aborted.
\begin{figure}[!htb]
	\centering
	{%
		% This file was created with tikzplotlib v0.10.1.
\begin{tikzpicture}[scale=0.7]

\definecolor{darkgray176}{RGB}{176,176,176}
\definecolor{darkorange25512714}{RGB}{255,127,14}
\definecolor{forestgreen4416044}{RGB}{44,160,44}
\definecolor{lightgray204}{RGB}{204,204,204}
\definecolor{steelblue31119180}{RGB}{31,119,180}

\begin{axis}[
legend cell align={left},
legend style={
  fill opacity=0.8,
  draw opacity=1,
  text opacity=1,
  at={(0.97,0.03)},
  anchor=south east,
  draw=lightgray204
},
tick align=outside,
tick pos=left,
title={objective function value},
x grid style={darkgray176},
xmin=0.9996535, xmax=1.0072765,
xtick style={color=black},
x tick label style={/pgf/number format/.cd, precision=3, /tikz/.cd,},
y grid style={darkgray176},
ymin=-0.05, ymax=1.05,
ytick style={color=black}
]
\addplot [semithick, steelblue31119180, const plot mark right]
table {%
1 0.666666666666667
1.00007 0.666666666666667
1.00014 0.666666666666667
1.00021 0.666666666666667
1.00028 0.666666666666667
1.00035 0.666666666666667
1.00042 0.666666666666667
1.00049 0.666666666666667
1.00056 0.666666666666667
1.00063 0.666666666666667
1.0007 0.666666666666667
1.00077 0.666666666666667
1.00084 0.666666666666667
1.00091 0.666666666666667
1.00098 0.666666666666667
1.00105 0.666666666666667
1.00112 0.666666666666667
1.00119 0.666666666666667
1.00126 0.666666666666667
1.00133 0.666666666666667
1.0014 0.666666666666667
1.00147 0.666666666666667
1.00154 0.666666666666667
1.00161 0.666666666666667
1.00168 0.666666666666667
1.00175 0.666666666666667
1.00182 0.666666666666667
1.00189 0.666666666666667
1.00196 0.666666666666667
1.00203 0.666666666666667
1.0021 0.666666666666667
1.00217 0.666666666666667
1.00224 0.666666666666667
1.00231 0.666666666666667
1.00238 0.666666666666667
1.00245 0.666666666666667
1.00252 0.666666666666667
1.00259 0.666666666666667
1.00266 0.666666666666667
1.00273 0.666666666666667
1.0028 0.666666666666667
1.00287 0.666666666666667
1.00294 0.666666666666667
1.00301 0.666666666666667
1.00308 0.666666666666667
1.00315 0.666666666666667
1.00322 0.666666666666667
1.00329 0.666666666666667
1.00336 0.666666666666667
1.00343 0.666666666666667
1.0035 0.666666666666667
1.00357 0.666666666666667
1.00364 0.666666666666667
1.00371 0.666666666666667
1.00378 0.666666666666667
1.00385 0.666666666666667
1.00392 0.666666666666667
1.00399 0.666666666666667
1.00406 0.666666666666667
1.00413 0.666666666666667
1.0042 0.666666666666667
1.00427 0.666666666666667
1.00434 0.666666666666667
1.00441 0.666666666666667
1.00448 0.666666666666667
1.00455 0.666666666666667
1.00462 0.7
1.00469 0.7
1.00476 0.7
1.00483 0.7
1.0049 0.7
1.00497 0.7
1.00504 0.7
1.00511 0.7
1.00518 0.7
1.00525 0.7
1.00532 0.7
1.00539 0.7
1.00546 0.7
1.00553 0.7
1.0056 0.766666666666667
1.00567 0.766666666666667
1.00574 0.766666666666667
1.00581 0.8
1.00588 0.833333333333333
1.00595 0.833333333333333
1.00602 0.833333333333333
1.00609 0.866666666666667
1.00616 0.9
1.00623 0.933333333333333
1.0063 0.933333333333333
1.00637 0.933333333333333
1.00644 0.933333333333333
1.00651 0.933333333333333
1.00658 0.933333333333333
1.00665 0.933333333333333
1.00672 0.933333333333333
1.00679 0.933333333333333
1.00686 0.933333333333333
1.00693 0.966666666666667
};
\addlegendentry{ALM}
\addplot [semithick, darkorange25512714, const plot mark right]
table {%
1 0
1.00007 0
1.00014 0
1.00021 0
1.00028 0.0333333333333333
1.00035 0.133333333333333
1.00042 0.266666666666667
1.00049 0.3
1.00056 0.366666666666667
1.00063 0.4
1.0007 0.4
1.00077 0.466666666666667
1.00084 0.5
1.00091 0.6
1.00098 0.633333333333333
1.00105 0.633333333333333
1.00112 0.766666666666667
1.00119 0.833333333333333
1.00126 0.833333333333333
1.00133 0.833333333333333
1.0014 0.9
1.00147 0.9
1.00154 0.9
1.00161 0.9
1.00168 0.9
1.00175 0.966666666666667
1.00182 1
1.00189 1
1.00196 1
1.00203 1
1.0021 1
1.00217 1
1.00224 1
1.00231 1
1.00238 1
1.00245 1
1.00252 1
1.00259 1
1.00266 1
1.00273 1
1.0028 1
1.00287 1
1.00294 1
1.00301 1
1.00308 1
1.00315 1
1.00322 1
1.00329 1
1.00336 1
1.00343 1
1.0035 1
1.00357 1
1.00364 1
1.00371 1
1.00378 1
1.00385 1
1.00392 1
1.00399 1
1.00406 1
1.00413 1
1.0042 1
1.00427 1
1.00434 1
1.00441 1
1.00448 1
1.00455 1
1.00462 1
1.00469 1
1.00476 1
1.00483 1
1.0049 1
1.00497 1
1.00504 1
1.00511 1
1.00518 1
1.00525 1
1.00532 1
1.00539 1
1.00546 1
1.00553 1
1.0056 1
1.00567 1
1.00574 1
1.00581 1
1.00588 1
1.00595 1
1.00602 1
1.00609 1
1.00616 1
1.00623 1
1.0063 1
1.00637 1
1.00644 1
1.00651 1
1.00658 1
1.00665 1
1.00672 1
1.00679 1
1.00686 1
1.00693 1
};
\addlegendentry{ALMr}
\addplot [semithick, forestgreen4416044, const plot mark right]
table {%
1 0.333333333333333
1.00007 0.333333333333333
1.00014 0.366666666666667
1.00021 0.366666666666667
1.00028 0.4
1.00035 0.6
1.00042 0.633333333333333
1.00049 0.733333333333333
1.00056 0.8
1.00063 0.866666666666667
1.0007 0.866666666666667
1.00077 0.9
1.00084 0.9
1.00091 0.933333333333333
1.00098 1
1.00105 1
1.00112 1
1.00119 1
1.00126 1
1.00133 1
1.0014 1
1.00147 1
1.00154 1
1.00161 1
1.00168 1
1.00175 1
1.00182 1
1.00189 1
1.00196 1
1.00203 1
1.0021 1
1.00217 1
1.00224 1
1.00231 1
1.00238 1
1.00245 1
1.00252 1
1.00259 1
1.00266 1
1.00273 1
1.0028 1
1.00287 1
1.00294 1
1.00301 1
1.00308 1
1.00315 1
1.00322 1
1.00329 1
1.00336 1
1.00343 1
1.0035 1
1.00357 1
1.00364 1
1.00371 1
1.00378 1
1.00385 1
1.00392 1
1.00399 1
1.00406 1
1.00413 1
1.0042 1
1.00427 1
1.00434 1
1.00441 1
1.00448 1
1.00455 1
1.00462 1
1.00469 1
1.00476 1
1.00483 1
1.0049 1
1.00497 1
1.00504 1
1.00511 1
1.00518 1
1.00525 1
1.00532 1
1.00539 1
1.00546 1
1.00553 1
1.0056 1
1.00567 1
1.00574 1
1.00581 1
1.00588 1
1.00595 1
1.00602 1
1.00609 1
1.00616 1
1.00623 1
1.0063 1
1.00637 1
1.00644 1
1.00651 1
1.00658 1
1.00665 1
1.00672 1
1.00679 1
1.00686 1
1.00693 1
};
\addlegendentry{sALMr}
\end{axis}

\end{tikzpicture}
		\hspace{1cm}
		% This file was created with tikzplotlib v0.10.1.
\begin{tikzpicture}[scale=0.7]

\definecolor{darkgray176}{RGB}{176,176,176}
\definecolor{darkorange25512714}{RGB}{255,127,14}
\definecolor{forestgreen4416044}{RGB}{44,160,44}
\definecolor{lightgray204}{RGB}{204,204,204}
\definecolor{steelblue31119180}{RGB}{31,119,180}

\begin{axis}[
legend cell align={left},
legend style={
  fill opacity=0.8,
  draw opacity=1,
  text opacity=1,
  at={(0.97,0.03)},
  anchor=south east,
  draw=lightgray204
},
tick align=outside,
tick pos=left,
title={number of outer iterations},
x grid style={darkgray176},
xmin=0.9802, xmax=1.4158,
xtick style={color=black},
y grid style={darkgray176},
ymin=-0.05, ymax=1.05,
ytick style={color=black}
]
\addplot [semithick, steelblue31119180, const plot mark right]
table {%
1 0.433333333333333
1.004 0.433333333333333
1.008 0.433333333333333
1.012 0.433333333333333
1.016 0.433333333333333
1.02 0.433333333333333
1.024 0.433333333333333
1.028 0.433333333333333
1.032 0.433333333333333
1.036 0.5
1.04 0.5
1.044 0.5
1.048 0.5
1.052 0.5
1.056 0.5
1.06 0.5
1.064 0.5
1.068 0.5
1.072 0.566666666666667
1.076 0.6
1.08 0.6
1.084 0.6
1.088 0.6
1.092 0.6
1.096 0.6
1.1 0.6
1.104 0.6
1.108 0.6
1.112 0.633333333333333
1.116 0.7
1.12 0.733333333333333
1.124 0.733333333333333
1.128 0.733333333333333
1.132 0.733333333333333
1.136 0.733333333333333
1.14 0.733333333333333
1.144 0.733333333333333
1.148 0.733333333333333
1.152 0.733333333333333
1.156 0.833333333333333
1.16 0.866666666666667
1.164 0.866666666666667
1.168 0.866666666666667
1.172 0.866666666666667
1.176 0.866666666666667
1.18 0.866666666666667
1.184 0.866666666666667
1.188 0.866666666666667
1.192 0.866666666666667
1.196 0.866666666666667
1.2 0.966666666666667
1.204 0.966666666666667
1.208 0.966666666666667
1.212 0.966666666666667
1.216 0.966666666666667
1.22 0.966666666666667
1.224 0.966666666666667
1.228 0.966666666666667
1.232 0.966666666666667
1.236 0.966666666666667
1.24 0.966666666666667
1.244 0.966666666666667
1.248 0.966666666666667
1.252 0.966666666666667
1.256 0.966666666666667
1.26 0.966666666666667
1.264 0.966666666666667
1.268 0.966666666666667
1.272 0.966666666666667
1.276 0.966666666666667
1.28 0.966666666666667
1.284 0.966666666666667
1.288 0.966666666666667
1.292 0.966666666666667
1.296 0.966666666666667
1.3 0.966666666666667
1.304 0.966666666666667
1.308 1
1.312 1
1.316 1
1.32 1
1.324 1
1.328 1
1.332 1
1.336 1
1.34 1
1.344 1
1.348 1
1.352 1
1.356 1
1.36 1
1.364 1
1.368 1
1.372 1
1.376 1
1.38 1
1.384 1
1.388 1
1.392 1
1.396 1
};
\addlegendentry{ALM}
\addplot [semithick, darkorange25512714, const plot mark right]
table {%
1 0.6
1.004 0.6
1.008 0.6
1.012 0.6
1.016 0.6
1.02 0.6
1.024 0.6
1.028 0.6
1.032 0.6
1.036 0.633333333333333
1.04 0.633333333333333
1.044 0.666666666666667
1.048 0.666666666666667
1.052 0.666666666666667
1.056 0.666666666666667
1.06 0.666666666666667
1.064 0.666666666666667
1.068 0.666666666666667
1.072 0.666666666666667
1.076 0.666666666666667
1.08 0.666666666666667
1.084 0.7
1.088 0.7
1.092 0.7
1.096 0.7
1.1 0.7
1.104 0.7
1.108 0.7
1.112 0.733333333333333
1.116 0.733333333333333
1.12 0.733333333333333
1.124 0.733333333333333
1.128 0.733333333333333
1.132 0.733333333333333
1.136 0.733333333333333
1.14 0.733333333333333
1.144 0.733333333333333
1.148 0.733333333333333
1.152 0.733333333333333
1.156 0.733333333333333
1.16 0.766666666666667
1.164 0.766666666666667
1.168 0.8
1.172 0.8
1.176 0.8
1.18 0.8
1.184 0.8
1.188 0.8
1.192 0.8
1.196 0.8
1.2 0.8
1.204 0.8
1.208 0.8
1.212 0.8
1.216 0.8
1.22 0.8
1.224 0.8
1.228 0.8
1.232 0.8
1.236 0.8
1.24 0.8
1.244 0.8
1.248 0.8
1.252 0.8
1.256 0.8
1.26 0.8
1.264 0.8
1.268 0.8
1.272 0.8
1.276 0.866666666666667
1.28 0.866666666666667
1.284 0.866666666666667
1.288 0.866666666666667
1.292 0.866666666666667
1.296 0.866666666666667
1.3 0.866666666666667
1.304 0.866666666666667
1.308 0.866666666666667
1.312 0.866666666666667
1.316 0.866666666666667
1.32 0.866666666666667
1.324 0.866666666666667
1.328 0.866666666666667
1.332 0.866666666666667
1.336 0.866666666666667
1.34 0.866666666666667
1.344 0.866666666666667
1.348 0.866666666666667
1.352 0.866666666666667
1.356 0.866666666666667
1.36 0.866666666666667
1.364 0.9
1.368 0.9
1.372 0.9
1.376 0.9
1.38 0.9
1.384 0.933333333333333
1.388 0.933333333333333
1.392 0.933333333333333
1.396 0.933333333333333
};
\addlegendentry{ALMr}
\addplot [semithick, forestgreen4416044, const plot mark right]
table {%
1 0
1.004 0
1.008 0
1.012 0
1.016 0
1.02 0
1.024 0
1.028 0
1.032 0
1.036 0.0666666666666667
1.04 0.0666666666666667
1.044 0.0666666666666667
1.048 0.0666666666666667
1.052 0.0666666666666667
1.056 0.0666666666666667
1.06 0.0666666666666667
1.064 0.0666666666666667
1.068 0.0666666666666667
1.072 0.166666666666667
1.076 0.166666666666667
1.08 0.166666666666667
1.084 0.166666666666667
1.088 0.166666666666667
1.092 0.166666666666667
1.096 0.166666666666667
1.1 0.166666666666667
1.104 0.166666666666667
1.108 0.166666666666667
1.112 0.3
1.116 0.3
1.12 0.3
1.124 0.3
1.128 0.3
1.132 0.3
1.136 0.3
1.14 0.3
1.144 0.3
1.148 0.3
1.152 0.3
1.156 0.466666666666667
1.16 0.466666666666667
1.164 0.466666666666667
1.168 0.466666666666667
1.172 0.466666666666667
1.176 0.466666666666667
1.18 0.466666666666667
1.184 0.466666666666667
1.188 0.466666666666667
1.192 0.466666666666667
1.196 0.466666666666667
1.2 0.666666666666667
1.204 0.666666666666667
1.208 0.666666666666667
1.212 0.666666666666667
1.216 0.666666666666667
1.22 0.666666666666667
1.224 0.666666666666667
1.228 0.666666666666667
1.232 0.666666666666667
1.236 0.666666666666667
1.24 0.666666666666667
1.244 0.666666666666667
1.248 0.666666666666667
1.252 0.766666666666667
1.256 0.766666666666667
1.26 0.766666666666667
1.264 0.766666666666667
1.268 0.766666666666667
1.272 0.766666666666667
1.276 0.766666666666667
1.28 0.766666666666667
1.284 0.766666666666667
1.288 0.766666666666667
1.292 0.766666666666667
1.296 0.766666666666667
1.3 0.766666666666667
1.304 0.766666666666667
1.308 0.8
1.312 0.8
1.316 0.8
1.32 0.8
1.324 0.8
1.328 0.8
1.332 0.8
1.336 0.8
1.34 0.8
1.344 0.8
1.348 0.8
1.352 0.8
1.356 0.8
1.36 0.8
1.364 0.9
1.368 0.9
1.372 0.9
1.376 0.9
1.38 0.9
1.384 0.9
1.388 0.9
1.392 0.9
1.396 0.9
};
\addlegendentry{sALMr}
\end{axis}

\end{tikzpicture}
	}%
	\\[.3cm]
	{%
		\input{value_of_error_measure.tex}
		\hspace{1cm}
		\input{value_of_penalty_parameter.tex}
	}%
	\caption{%
		Performance profiles for 
		the objective function value, 
		the number of outer augmented Lagrangian iterations,
		the value of the error measure,
		and the value of the penalty parameter
		(from top left to bottom right).
		}%
	\label{fig:basic_PPs}
\end{figure}
We observe that \ALM\ computes the best objective function values for approximately 67\% of the test instances,
but the scale of the $t$-axis underlines that \ALMR\ and \SALMR\ do not significantly fall behind \ALM\
in this regard. Interestingly, \ALMR\ and \SALMR\ produce a point whose objective
value is just 1\textperthousand\ and 2\textperthousand\ worse than the best computed function value for
each test instance, respectively, a quality which \ALM\ does not achieve.
Regarding the total number of outer augmented Lagrangian iterations, we observe that \ALM\ and \ALMR\
perform best in about 40\% and 60\% of all test runs, respectively, 
while \SALMR\ cannot match the other two methods in this regard.
The performance profile also indicates that in about 90\% of all test runs, the total number of outer iterations
for any of the algorithms under consideration
is at most $1.5$ times as high as the smallest number of exploited outer iterations.
Regarding the value of the error measure, there is a significant difference between \ALM\ as well as \ALMR,
which perform similarly good with slight advantages for \ALMR, and \SALMR\ that is outperformed
by the other two methods. Exemplary, the situation where the error measure is allowed to be 20 times higher
than the best achieved value of the error measure is only covered in about 60\% of all test runs for \SALMR.
The final value of the penalty parameter for \ALM\ is smaller than the one for \ALMR\ for about
68\% of all test instances. However, while, for its termination, \ALMR\ reliably needs a penalty parameter which is at most
$\gamma^6=4096$ times larger than the smallest necessary penalty parameter, i.e., at most $6$ more enlargements
of the penalty parameter in \cref{item:update_penalty_ALM} of \cref{alg:ALM} are necessary, 
a similar bound for \ALM\ cannot be distilled from the associated performance profile.

Second, let us inspect the feasibility properties of the outputs of all these methods.
In \cref{fig:PP_percentage_violated_constraints}, the percentage of violated constraints
(w.r.t.\ the model \eqref{eq:VPD}) is documented.
\begin{figure}[!htb]
	\centering
	{%
		\input{percentage_of_violated_constraints.tex}
	}%
	\caption{%
		Performance profile 
		for the percentage of violated constraints.
		}%
	\label{fig:PP_percentage_violated_constraints}
\end{figure}
Here, we observe that the smallest percentage of violated constraints
is achieved by \ALMR\ in about 68\% of all test runs. However, \ALM\ falls behind the smallest
percentage of violated constraints only by factor $1.003$ in all test runs, and, anyhow, is in the lead
for the other 32\% of the test instances. Again, \SALMR\ clearly cannot compete with \ALM\ and \ALMR.
Still, it falls short the smallest percentage of violated constraints only by factor $1.01$ which is
still an acceptable deviation. Note, however, that the percentage of violated constraints does not
document how bad the violated constraints are missed.
Therefore, let us inspect \cref{fig:PPs_average_relative_constraint_violation} which reports the
maximum and average relative (affine) constraint violation.
\begin{figure}[!htb]
	\centering
	{%
		\input{maximum_relative_constraint_violation.tex}
		\hspace{1cm}
		\input{maximum_relative_affine_constraint_violation.tex}
	}%
	\\[.3cm]
	{%
		% This file was created with tikzplotlib v0.10.1.
\begin{tikzpicture}[scale=0.7]

\definecolor{darkgray176}{RGB}{176,176,176}
\definecolor{darkorange25512714}{RGB}{255,127,14}
\definecolor{forestgreen4416044}{RGB}{44,160,44}
\definecolor{lightgray204}{RGB}{204,204,204}
\definecolor{steelblue31119180}{RGB}{31,119,180}

\begin{axis}[
legend cell align={left},
legend style={
  fill opacity=0.8,
  draw opacity=1,
  text opacity=1,
  at={(0.97,0.03)},
  anchor=south east,
  draw=lightgray204
},
tick align=outside,
tick pos=left,
title={average relative constraint violation},
x grid style={darkgray176},
xmin=0.9996535, xmax=1.0072765,
xtick style={color=black},
x tick label style={/pgf/number format/.cd, precision=3, /tikz/.cd,},
y grid style={darkgray176},
ymin=-0.05, ymax=1.05,
ytick style={color=black}
]
\addplot [semithick, steelblue31119180, const plot mark right]
table {%
1 0.333333333333333
1.00007 0.533333333333333
1.00014 0.866666666666667
1.00021 0.966666666666667
1.00028 1
1.00035 1
1.00042 1
1.00049 1
1.00056 1
1.00063 1
1.0007 1
1.00077 1
1.00084 1
1.00091 1
1.00098 1
1.00105 1
1.00112 1
1.00119 1
1.00126 1
1.00133 1
1.0014 1
1.00147 1
1.00154 1
1.00161 1
1.00168 1
1.00175 1
1.00182 1
1.00189 1
1.00196 1
1.00203 1
1.0021 1
1.00217 1
1.00224 1
1.00231 1
1.00238 1
1.00245 1
1.00252 1
1.00259 1
1.00266 1
1.00273 1
1.0028 1
1.00287 1
1.00294 1
1.00301 1
1.00308 1
1.00315 1
1.00322 1
1.00329 1
1.00336 1
1.00343 1
1.0035 1
1.00357 1
1.00364 1
1.00371 1
1.00378 1
1.00385 1
1.00392 1
1.00399 1
1.00406 1
1.00413 1
1.0042 1
1.00427 1
1.00434 1
1.00441 1
1.00448 1
1.00455 1
1.00462 1
1.00469 1
1.00476 1
1.00483 1
1.0049 1
1.00497 1
1.00504 1
1.00511 1
1.00518 1
1.00525 1
1.00532 1
1.00539 1
1.00546 1
1.00553 1
1.0056 1
1.00567 1
1.00574 1
1.00581 1
1.00588 1
1.00595 1
1.00602 1
1.00609 1
1.00616 1
1.00623 1
1.0063 1
1.00637 1
1.00644 1
1.00651 1
1.00658 1
1.00665 1
1.00672 1
1.00679 1
1.00686 1
1.00693 1
};
\addlegendentry{ALM}
\addplot [semithick, darkorange25512714, const plot mark right]
table {%
1 0.666666666666667
1.00007 0.666666666666667
1.00014 0.666666666666667
1.00021 0.666666666666667
1.00028 0.666666666666667
1.00035 0.666666666666667
1.00042 0.666666666666667
1.00049 0.666666666666667
1.00056 0.666666666666667
1.00063 0.7
1.0007 0.7
1.00077 0.7
1.00084 0.733333333333333
1.00091 0.766666666666667
1.00098 0.766666666666667
1.00105 0.766666666666667
1.00112 0.766666666666667
1.00119 0.8
1.00126 0.8
1.00133 0.8
1.0014 0.8
1.00147 0.833333333333333
1.00154 0.833333333333333
1.00161 0.9
1.00168 0.9
1.00175 0.9
1.00182 0.9
1.00189 0.933333333333333
1.00196 0.933333333333333
1.00203 0.933333333333333
1.0021 0.933333333333333
1.00217 0.933333333333333
1.00224 0.966666666666667
1.00231 0.966666666666667
1.00238 0.966666666666667
1.00245 0.966666666666667
1.00252 0.966666666666667
1.00259 0.966666666666667
1.00266 1
1.00273 1
1.0028 1
1.00287 1
1.00294 1
1.00301 1
1.00308 1
1.00315 1
1.00322 1
1.00329 1
1.00336 1
1.00343 1
1.0035 1
1.00357 1
1.00364 1
1.00371 1
1.00378 1
1.00385 1
1.00392 1
1.00399 1
1.00406 1
1.00413 1
1.0042 1
1.00427 1
1.00434 1
1.00441 1
1.00448 1
1.00455 1
1.00462 1
1.00469 1
1.00476 1
1.00483 1
1.0049 1
1.00497 1
1.00504 1
1.00511 1
1.00518 1
1.00525 1
1.00532 1
1.00539 1
1.00546 1
1.00553 1
1.0056 1
1.00567 1
1.00574 1
1.00581 1
1.00588 1
1.00595 1
1.00602 1
1.00609 1
1.00616 1
1.00623 1
1.0063 1
1.00637 1
1.00644 1
1.00651 1
1.00658 1
1.00665 1
1.00672 1
1.00679 1
1.00686 1
1.00693 1
};
\addlegendentry{ALMr}
\addplot [semithick, forestgreen4416044, const plot mark right]
table {%
1 0
1.00007 0
1.00014 0.0666666666666667
1.00021 0.3
1.00028 0.566666666666667
1.00035 0.666666666666667
1.00042 0.666666666666667
1.00049 0.666666666666667
1.00056 0.666666666666667
1.00063 0.666666666666667
1.0007 0.666666666666667
1.00077 0.666666666666667
1.00084 0.666666666666667
1.00091 0.666666666666667
1.00098 0.666666666666667
1.00105 0.666666666666667
1.00112 0.666666666666667
1.00119 0.666666666666667
1.00126 0.666666666666667
1.00133 0.666666666666667
1.0014 0.666666666666667
1.00147 0.666666666666667
1.00154 0.666666666666667
1.00161 0.666666666666667
1.00168 0.666666666666667
1.00175 0.7
1.00182 0.7
1.00189 0.7
1.00196 0.7
1.00203 0.7
1.0021 0.7
1.00217 0.7
1.00224 0.7
1.00231 0.7
1.00238 0.7
1.00245 0.733333333333333
1.00252 0.733333333333333
1.00259 0.733333333333333
1.00266 0.733333333333333
1.00273 0.733333333333333
1.0028 0.733333333333333
1.00287 0.733333333333333
1.00294 0.733333333333333
1.00301 0.733333333333333
1.00308 0.733333333333333
1.00315 0.733333333333333
1.00322 0.733333333333333
1.00329 0.733333333333333
1.00336 0.733333333333333
1.00343 0.733333333333333
1.0035 0.733333333333333
1.00357 0.766666666666667
1.00364 0.8
1.00371 0.8
1.00378 0.8
1.00385 0.8
1.00392 0.8
1.00399 0.8
1.00406 0.8
1.00413 0.8
1.0042 0.8
1.00427 0.8
1.00434 0.8
1.00441 0.833333333333333
1.00448 0.833333333333333
1.00455 0.866666666666667
1.00462 0.9
1.00469 0.9
1.00476 0.9
1.00483 0.9
1.0049 0.933333333333333
1.00497 0.933333333333333
1.00504 0.933333333333333
1.00511 0.933333333333333
1.00518 0.933333333333333
1.00525 0.933333333333333
1.00532 0.933333333333333
1.00539 0.966666666666667
1.00546 0.966666666666667
1.00553 0.966666666666667
1.0056 0.966666666666667
1.00567 0.966666666666667
1.00574 0.966666666666667
1.00581 0.966666666666667
1.00588 0.966666666666667
1.00595 1
1.00602 1
1.00609 1
1.00616 1
1.00623 1
1.0063 1
1.00637 1
1.00644 1
1.00651 1
1.00658 1
1.00665 1
1.00672 1
1.00679 1
1.00686 1
1.00693 1
};
\addlegendentry{sALMr}
\end{axis}

\end{tikzpicture}
		\hspace{1cm}
		% This file was created with tikzplotlib v0.10.1.
\begin{tikzpicture}[scale=0.7]

\definecolor{darkgray176}{RGB}{176,176,176}
\definecolor{darkorange25512714}{RGB}{255,127,14}
\definecolor{forestgreen4416044}{RGB}{44,160,44}
\definecolor{lightgray204}{RGB}{204,204,204}
\definecolor{steelblue31119180}{RGB}{31,119,180}

\begin{axis}[
legend cell align={left},
legend style={
  fill opacity=0.8,
  draw opacity=1,
  text opacity=1,
  at={(0.97,0.03)},
  anchor=south east,
  draw=lightgray204
},
tick align=outside,
tick pos=left,
title={average relative affine constraint violation},
x grid style={darkgray176},
xmin=0.99998515, xmax=1.00031185,
xtick style={color=black},
xtick distance={0.0001},
x tick label style={/pgf/number format/.cd, precision=5, /tikz/.cd,},
y grid style={darkgray176},
ymin=-0.05, ymax=1.05,
ytick style={color=black}
]
\addplot [semithick, steelblue31119180, const plot mark right]
table {%
1 0.333333333333333
1.000003 0.333333333333333
1.000006 0.333333333333333
1.000009 0.333333333333333
1.000012 0.333333333333333
1.000015 0.333333333333333
1.000018 0.333333333333333
1.000021 0.366666666666667
1.000024 0.4
1.000027 0.433333333333333
1.00003 0.466666666666667
1.000033 0.566666666666667
1.000036 0.6
1.000039 0.633333333333333
1.000042 0.666666666666667
1.000045 0.7
1.000048 0.733333333333333
1.000051 0.733333333333333
1.000054 0.766666666666667
1.000057 0.866666666666667
1.00006 0.866666666666667
1.000063 0.866666666666667
1.000066 0.866666666666667
1.000069 0.9
1.000072 0.9
1.000075 0.9
1.000078 0.9
1.000081 0.9
1.000084 0.933333333333333
1.000087 0.933333333333333
1.00009 0.933333333333333
1.000093 0.966666666666667
1.000096 0.966666666666667
1.000099 0.966666666666667
1.000102 0.966666666666667
1.000105 0.966666666666667
1.000108 1
1.000111 1
1.000114 1
1.000117 1
1.00012 1
1.000123 1
1.000126 1
1.000129 1
1.000132 1
1.000135 1
1.000138 1
1.000141 1
1.000144 1
1.000147 1
1.00015 1
1.000153 1
1.000156 1
1.000159 1
1.000162 1
1.000165 1
1.000168 1
1.000171 1
1.000174 1
1.000177 1
1.00018 1
1.000183 1
1.000186 1
1.000189 1
1.000192 1
1.000195 1
1.000198 1
1.000201 1
1.000204 1
1.000207 1
1.00021 1
1.000213 1
1.000216 1
1.000219 1
1.000222 1
1.000225 1
1.000228 1
1.000231 1
1.000234 1
1.000237 1
1.00024 1
1.000243 1
1.000246 1
1.000249 1
1.000252 1
1.000255 1
1.000258 1
1.000261 1
1.000264 1
1.000267 1
1.00027 1
1.000273 1
1.000276 1
1.000279 1
1.000282 1
1.000285 1
1.000288 1
1.000291 1
1.000294 1
1.000297 1
};
\addlegendentry{ALM}
\addplot [semithick, darkorange25512714, const plot mark right]
table {%
1 0.666666666666667
1.000003 1
1.000006 1
1.000009 1
1.000012 1
1.000015 1
1.000018 1
1.000021 1
1.000024 1
1.000027 1
1.00003 1
1.000033 1
1.000036 1
1.000039 1
1.000042 1
1.000045 1
1.000048 1
1.000051 1
1.000054 1
1.000057 1
1.00006 1
1.000063 1
1.000066 1
1.000069 1
1.000072 1
1.000075 1
1.000078 1
1.000081 1
1.000084 1
1.000087 1
1.00009 1
1.000093 1
1.000096 1
1.000099 1
1.000102 1
1.000105 1
1.000108 1
1.000111 1
1.000114 1
1.000117 1
1.00012 1
1.000123 1
1.000126 1
1.000129 1
1.000132 1
1.000135 1
1.000138 1
1.000141 1
1.000144 1
1.000147 1
1.00015 1
1.000153 1
1.000156 1
1.000159 1
1.000162 1
1.000165 1
1.000168 1
1.000171 1
1.000174 1
1.000177 1
1.00018 1
1.000183 1
1.000186 1
1.000189 1
1.000192 1
1.000195 1
1.000198 1
1.000201 1
1.000204 1
1.000207 1
1.00021 1
1.000213 1
1.000216 1
1.000219 1
1.000222 1
1.000225 1
1.000228 1
1.000231 1
1.000234 1
1.000237 1
1.00024 1
1.000243 1
1.000246 1
1.000249 1
1.000252 1
1.000255 1
1.000258 1
1.000261 1
1.000264 1
1.000267 1
1.00027 1
1.000273 1
1.000276 1
1.000279 1
1.000282 1
1.000285 1
1.000288 1
1.000291 1
1.000294 1
1.000297 1
};
\addlegendentry{ALMr}
\addplot [semithick, forestgreen4416044, const plot mark right]
table {%
1 0
1.000003 0.133333333333333
1.000006 0.333333333333333
1.000009 0.333333333333333
1.000012 0.333333333333333
1.000015 0.333333333333333
1.000018 0.333333333333333
1.000021 0.333333333333333
1.000024 0.333333333333333
1.000027 0.333333333333333
1.00003 0.333333333333333
1.000033 0.333333333333333
1.000036 0.333333333333333
1.000039 0.333333333333333
1.000042 0.333333333333333
1.000045 0.333333333333333
1.000048 0.333333333333333
1.000051 0.333333333333333
1.000054 0.333333333333333
1.000057 0.333333333333333
1.00006 0.333333333333333
1.000063 0.333333333333333
1.000066 0.333333333333333
1.000069 0.333333333333333
1.000072 0.333333333333333
1.000075 0.333333333333333
1.000078 0.333333333333333
1.000081 0.366666666666667
1.000084 0.366666666666667
1.000087 0.366666666666667
1.00009 0.4
1.000093 0.433333333333333
1.000096 0.466666666666667
1.000099 0.5
1.000102 0.533333333333333
1.000105 0.6
1.000108 0.6
1.000111 0.6
1.000114 0.6
1.000117 0.633333333333333
1.00012 0.633333333333333
1.000123 0.633333333333333
1.000126 0.633333333333333
1.000129 0.633333333333333
1.000132 0.633333333333333
1.000135 0.633333333333333
1.000138 0.666666666666667
1.000141 0.666666666666667
1.000144 0.666666666666667
1.000147 0.666666666666667
1.00015 0.666666666666667
1.000153 0.666666666666667
1.000156 0.666666666666667
1.000159 0.666666666666667
1.000162 0.666666666666667
1.000165 0.666666666666667
1.000168 0.666666666666667
1.000171 0.666666666666667
1.000174 0.666666666666667
1.000177 0.666666666666667
1.00018 0.666666666666667
1.000183 0.666666666666667
1.000186 0.666666666666667
1.000189 0.666666666666667
1.000192 0.666666666666667
1.000195 0.666666666666667
1.000198 0.666666666666667
1.000201 0.666666666666667
1.000204 0.666666666666667
1.000207 0.666666666666667
1.00021 0.666666666666667
1.000213 0.666666666666667
1.000216 0.733333333333333
1.000219 0.766666666666667
1.000222 0.8
1.000225 0.833333333333333
1.000228 0.833333333333333
1.000231 0.866666666666667
1.000234 0.866666666666667
1.000237 0.9
1.00024 0.9
1.000243 0.933333333333333
1.000246 1
1.000249 1
1.000252 1
1.000255 1
1.000258 1
1.000261 1
1.000264 1
1.000267 1
1.00027 1
1.000273 1
1.000276 1
1.000279 1
1.000282 1
1.000285 1
1.000288 1
1.000291 1
1.000294 1
1.000297 1
};
\addlegendentry{sALMr}
\end{axis}

\end{tikzpicture}
	}%
	\caption{%
		Performance profiles for 
		the maximum relative constraint violation, 
		the maximum relative affine constraint violation,
		the average relative constraint violation,
		and the average relative affine constraint violation
		(from top left to bottom right).
		}%
	\label{fig:PPs_average_relative_constraint_violation}
\end{figure}
Here, we observe that regarding the maximum and average relative constraint violation
(w.r.t.\ the constraints of \eqref{eq:VPD}), all methods perform in a good way
with obvious advantages for \ALM\ and \ALMR\ in the settings of the average and
maximum constraint violation, respectively.
This is good news, as all three methods are, at their core, constructed to solve
the model problem \eqref{eq:VPD} which comes along with the handling of the
associated huge number of constraints.
Let us now inspect the performance profiles for the maximum and average relative affine constraint violation
(w.r.t.\ the constraints of \eqref{eq:VPDaff}).
One would await that \ALMR\ as well as \SALMR\ outperform \ALM\ in this regard as the former
two algorithms directly work with the model formulation \eqref{eq:VPDaff} 
while the latter algorithm does not. This intuition is supported by our results only in parts.
Indeed, we observe that \ALM\ and \SALMR\ both cannot challenge \ALMR,
and that \ALM\ performs significantly better than \SALMR.

Let us present an interim summary of our numerical experiments.
While \ALM\ and \ALMR\ both perform well on the set of benchmark problems,
\SALMR\ falls short both of these methods regarding the number of outer iterations,
the final value of the error measure, and the feasibility aspects of the produced output.
Furthermore, it has to be mentioned that, due to the fact that we need to run 10000
iterations of NADAM in order to solve the first augmented Lagrangian subproblem up to
a reasonable precision to observe any convergence of \SALMR\ later on,
see \cref{sec:denoising_SGD}, the latter is significantly slower than \ALM\ and \ALMR.
We, thus, cannot attest \SALMR\ a satisfying performance.

Next, we present some averaged numbers to analyze the behavior of the three methods
on each of the three test images individually. 
Let us start with a documentation of the results associated with the test image ``brain''
which can be found in \cref{tab:brain}.
\begin{table}[htp]
\centering
\begin{tabular}{llll}
\toprule
performance measure & \ALM\ & \ALMR\ & \SALMR\
\\
\midrule
objective function value 				& 2.4729$\cdot 10^{8}$		& 2.4592$\cdot 10^{8}$ 		& 2.4583$\cdot 10^{8}$ \\
number of outer iterations 				& 29.3 						& 26.5 						& 30.0 \\
value of error measure 					& 1.4998$\cdot 10^{1}$ 		& 7.6307$\cdot 10^{-3}$ 	& 3.8705$\cdot 10^{-2}$ \\
value of penalty parameter 				& 3.4903$\cdot 10^{16}$		& 1.0486$\cdot 10^{7}$ 		& 4.0000$\cdot 10^{5}$ \\
percentage of violated constraints 		& 9.4996$\cdot 10^{-5}$ 	& 5.8363$\cdot 10^{-3}$ 	& 8.4727$\cdot 10^{-3}$ \\
max.\ rel.\ constraint violation 		& 3.7626$\cdot 10^{-3}$ 	& 9.6763$\cdot 10^{-1}$ 	& 1.7180 \\
max.\ rel.\ affine constraint violation & 2.2717$\cdot 10^{-5}$ 	& 7.6307$\cdot 10^{-3}$ 	& 3.8705$\cdot 10^{-2}$ \\
av.\ rel.\ constraint violation 		& 1.1616$\cdot 10^{-8}$ 	& 1.4712$\cdot 10^{-3}$ 	& 4.0849$\cdot 10^{-3}$ \\
av.\ rel.\ aff.\ constraint violation 	& 2.2768$\cdot 10^{-11}$ 	& 6.7136$\cdot 10^{-7}$ 	& 3.3754$\cdot 10^{-6}$ \\
\bottomrule
\end{tabular}
\caption{%
	Averaged numbers for the test image ``brain''.
}\label{tab:brain}
\end{table}
The performance profiles already indicated that all three methods compute points with
related objective function value, and this is underlined by \cref{tab:brain}.
Inspecting the total number of outer iterations reveals an advantage for \ALMR\ against
\ALM\ as well as \SALMR, and similar observations can be made for the final value of 
error measure and penalty parameter. 
Let us note that \SALMR\ exploits the predefined maximum total number of outer iterations 
in each run. 
It, thus, might be possible to achieve better results with this approach when allowing for
further iterations. This, however, costs additional time, and we already mentioned above
that \SALMR\ is the slowest algorithm we are considering here.
Regarding the feasibility aspects,
\ALM\ is clearly better than the other two methods.
It has to be mentioned that there is a strong relation between these assessments.
For the test image ``brain'', which possesses broad areas of uniform color, 
we observe that, on many subboxes,
the mean of the noisy image as well as the mean of the iterates are (very close to) zero. 
Numerically, this drives the value of the associated constraints in \eqref{eq:VPD}
and, thus, the error measure to $\infty$. In order to avoid this, \cref{alg:ALM}
enlarges the penalty parameter more often in this setting, compared to the one
where \eqref{eq:VPDaff} is tackled and this phenomenon cannot happen.
At the end of the day, \ALM\ then produces images which are much closer to
feasibility.

Let us now inspect the averaged numbers for the results associated with the
test images ``butterfly'' and ``cameraman'' in \cref{tab:butterfly,tab:cameraman},
respectively.
\begin{table}[htp]
\centering
\begin{tabular}{llll}
\toprule
performance measure & \ALM\ & \ALMR\ & \SALMR\
\\
\midrule
objective function value 				& 9.8967$\cdot 10^{8}$		& 9.9047$\cdot 10^{8}$ 		& 9.9007$\cdot 10^{8}$ \\
number of outer iterations 				& 27.9 						& 26.4						& 30.0 \\
value of error measure 					& 8.8339$\cdot 10^{-3}$ 	& 7.3673$\cdot 10^{-3}$ 	& 8.7777$\cdot 10^{-2}$ \\
value of penalty parameter 				& 1.2452$\cdot 10^{5}$ 		& 6.2076$\cdot 10^{7}$ 		& 4.0000$\cdot 10^{5}$ \\
percentage of violated constraints 		& 2.2192$\cdot 10^{-3}$ 	& 7.4833$\cdot 10^{-5}$ 	& 5.4026$\cdot 10^{-3}$ \\
max.\ rel.\ constraint violation 		& 9.6253$\cdot 10^{-1}$ 	& 1.8556$\cdot 10^{-2}$		& 3.8648$\cdot 10^{-1}$ \\
max.\ rel.\ affine constraint violation & 2.7286$\cdot 10^{-1}$		& 7.3673$\cdot 10^{-3}$		& 1.4300 \\
av.\ rel.\ constraint violation 		& 9.5834$\cdot 10^{-5}$ 	& 7.6614$\cdot 10^{-8}$ 	& 2.5851$\cdot 10^{-4}$ \\
av.\ rel.\ aff.\ constraint violation 	& 4.4251$\cdot 10^{-5}$ 	& 4.8633$\cdot 10^{-8}$ 	& 2.2870$\cdot 10^{-4}$ \\
\bottomrule
\end{tabular}
\caption{%
	Averaged numbers for the test image ``butterfly''.
}\label{tab:butterfly}
\end{table}

\begin{table}[htp]
\centering
\begin{tabular}{llll}
\toprule
performance measure & \ALM\ & \ALMR\ & \SALMR\
\\
\midrule
objective function value 				& 1.1606$\cdot 10^{9}$		& 1.1622$\cdot 10^{9}$ 		& 1.1613$\cdot 10^{9}$ \\
number of outer iterations 				& 23.8 						& 28.3						& 30.0 \\
value of error measure 					& 6.8892$\cdot 10^{-3}$ 	& 4.9436$\cdot 10^{-3}$ 	& 1.2357$\cdot 10^{-1}$ \\
value of penalty parameter 				& 3.0147$\cdot 10^{5}$ 		& 1.4764$\cdot 10^{8}$ 		& 4.0000$\cdot 10^{5}$ \\
percentage of violated constraints 		& 1.9836$\cdot 10^{-3}$ 	& 3.4649$\cdot 10^{-5}$ 	& 4.1128$\cdot 10^{-3}$ \\
max.\ rel.\ constraint violation 		& 7.7851$\cdot 10^{-1}$ 	& 1.2058$\cdot 10^{-2}$		& 3.4789$\cdot 10^{-1}$ \\
max.\ rel.\ affine constraint violation & 2.6348$\cdot 10^{-1}$		& 4.9436$\cdot 10^{-3}$		& 7.9586$\cdot 10^{-1}$ \\
av.\ rel.\ constraint violation 		& 1.1465$\cdot 10^{-4}$ 	& 5.1053$\cdot 10^{-8}$ 	& 1.6282$\cdot 10^{-4}$ \\
av.\ rel.\ aff.\ constraint violation 	& 5.1804$\cdot 10^{-5}$ 	& 2.5188$\cdot 10^{-8}$ 	& 1.0092$\cdot 10^{-4}$ \\
\bottomrule
\end{tabular}
\caption{%
	Averaged numbers for the test image ``cameraman''.
}\label{tab:cameraman}
\end{table}
Again, all three methods do not differ regarding the objective function value of the produced output.
Inspecting the error measure reveals no big differences between \ALM\ and \ALMR,
while \ALM\ now terminates with a significantly smaller penalty parameter than \ALMR\ does.
In contrast, the feasibility measures look far better for \ALMR\ than for \ALM.
Interestingly, for the test image ``cameraman'', \ALM\ terminates after a notably smaller number
of outer iterations than \ALMR.
As above, \SALMR\ exploits the predefined maximum total number of outer iterations in each run.
Regarding almost all aspects, \SALMR\ cannot match the other two methods.

Summing up these results, we observe that both \ALM\ and \ALMR\ are reasonable strategies 
for the implementation of the Poisson denoising approach for image recovery while \SALMR\ is not.
As depicted above, the individual performance of \ALM\ and \ALMR\ depends on the structure of
the underlying picture, so a general ranking between \ALM\ and \ALMR\ is not possible.
This is rather remarkable since \ALM\ is based on the nonlinear nonsmooth model problem
\eqref{eq:VPD} while \ALMR\ builds on the seemingly simpler convex quadratic problem \eqref{eq:VPDaff}.

\section{Concluding remarks}\label{sec:conclusions}

In this paper, we discussed the numerical solution of the so-called
variational Poisson denoising problem,
which is a statistically motivated model for image denoising that
possesses a huge number of constraints, with the aid of 
augmented Lagrangian methods.
To tackle the original model \eqref{eq:VPD}, 
whose constraints are modeled with the aid
of the extended real-valued Kullback--Leibler divergence
and are, thus, nonsmooth,
we suggested a rather general safeguarded augmented Lagrangian framework 
for fully nonsmooth problems in Banach spaces with finitely many inequality constraints, 
equality constraints within a Hilbert space setting, 
and additional abstract constraints,
where the inequality and equality constraints are augmented.
An associated derivative-free global convergence theory has been developed which applies 
in situations where the appearing subproblems can be solved to approximate global minimality, 
and the latter is likely to be possible in convex situations.
Our results generalize related findings in (partially) smooth settings,
see e.g.\ \cite[Section~4]{KanzowSteckWachsmuth2018} or \cite[Theorem~6.15]{KrugerMehlitz2022}.
For our analysis, we only relied on minimal requirements regarding semicontinuity properties
of all involved data functions as (generalized) differentiation played no role, and this makes
our results broadly applicable.

We also visualized that it is possible to equivalently reformulate the constraints
in \eqref{eq:VPD} as affine inequalities, but this procedure comes at the price of
doubling the total number of constraints. After discretization, this new model
\eqref{eq:VPDaff} is a convex quadratic optimization problem with a huge number of constraints,
and can be tackled with the standard safeguarded augmented Lagrangian framework.
Furthermore, as the minimizers of \eqref{eq:VPDaff} are stationary, it is possible
to solve the latter with the classical non-safeguarded version of the augmented Lagrangian method
where the penalty parameter undergoes no evolution.

These three algorithmic approaches have been implemented to denoise a benchmark collection
of noisy images obtained from the standard test images ``brain'', ``butterfly'', and
``cameraman'' by adding some random Poisson noise.
On the one hand, our results document that the safeguarded augmented Lagrangian method performs
comparably good when applied to \eqref{eq:VPD} or its reformulation \eqref{eq:VPDaff}.
On the other hand, the classical non-safeguarded augmented Lagrangian method 
with a constant penalty parameter falls clearly short of the results obtained 
by the other two methods and, thus, cannot be considered a reasonable approach to solve
the Poisson denoising model.

\subsection*{Acknowledgments}

The authors sincerely thank an anonymous reviewer of an earlier version of this paper
who suggested the investigation of the reformulation \eqref{eq:VPDaff}.
The research of Christian Kanzow and Gerd Wachsmuth 
was supported by the German Research Foundation (DFG) within the
priority program ``Non-smooth and Complementarity-based Distributed Parameter Systems: Simulation
and Hierarchical Optimization'' (SPP 1962) under grant numbers KA 1296/24-2 and WA 3636/4-2, respectively.

%%%%%% Bibliography
%%\bibliographystyle{plainnat}
%\bibliographystyle{habbrv}
%%\bibliographystyle{habbrv_fj}
%\bibliography{references}

\end{document}